\documentclass[]{article}

\usepackage{amssymb}
\usepackage{amsmath}

\usepackage{graphicx,subfigure}
\usepackage{amsmath, amsfonts, amsthm, amssymb, amstext, multicol, wrapfig, caption, , float, multirow, tabularx,bm,ragged2e,array,scrextend, verbatim,enumerate,mathtools,indentfirst}
\usepackage[colorlinks]{hyperref}
\hypersetup{
	citecolor = blue
}
\usepackage[dvipsnames]{xcolor}
\usepackage{authblk}

\theoremstyle{definition}
\usepackage[ruled]{algorithm2e}
\SetKwInput{KwInput}{Inputs}                
\SetKwInput{KwOutput}{Output} 
\SetKwComment{Comment}{$\triangleright$ }{}	

\newcommand{\id}{\mathrm d}
\newcommand{\vc}{\boldsymbol}
\newcommand{\bs}{\boldsymbol}

\newcommand{\pard}[2]{\frac{\partial #1}{\partial #2}}

\newcommand{\R}{\mathbb{R}}

\newcommand{\N}{\mathbb{N}}

\newtheorem{thm}{Theorem}
\newtheorem{rem}{Remark}
\newtheorem{lem}{Lemma}

\newtheorem{ass}{Assumption}

\newcommand{\longcomment}[1]{}

\newcommand{\singparam}{\theta}
\newcommand{\params}{\boldsymbol \singparam}
\newcommand{\numsens}{r}
\newcommand{\data}{y}
\newcommand{\SMSTikh}{\gamma}
\newcommand{\DATikh}{\tilde{\gamma}}

\newcommand{\ADEvel}{\vc v}
\newcommand{\ADEvelx}{v_1}
\newcommand{\ADEvelz}{v_2}
\newcommand{\ADEvar}{u} 
\newcommand{\ADEconst}{\kappa}

\newcommand{\hilspace}{\mathcal{H}}
\newcommand{\hilnorm}[1]{\|#1 \|_{\hilspace}}
\newcommand{\hilip}[2]{\left\langle#1, #2\right\rangle _{\hilspace}}

\newcommand{\spatdom}{\Omega}
\newcommand{\nummodes}{N}
\newcommand{\numcoll}{N_c}
\newcommand{\numparams}{p}
\newcommand{\numtimes}{m}

\newcommand{\jacop}{J_\numsens}

\begin{document}
	\title{Sequential data assimilation for PDEs using shape-morphing solutions\thanks{Accepted for publication in the Journal of Computational Physics. Corresponding author: farazmand@ncsu.edu}}
	\author[]{Zachary T. Hilliard} 
	\author[]{Mohammad Farazmand}
	\affil[]{Department of Mathematics, North Carolina State University, 2311 Stinson Drive, Raleigh, NC 27695, USA} 
	\renewcommand\Affilfont{\itshape\small}
	\date{}
	
\maketitle

\begin{abstract}
Shape-morphing solutions (also known as evolutional deep neural networks, reduced-order nonlinear solutions, and neural Galerkin schemes) are a new class of methods for approximating the solution of time-dependent partial differential equations (PDEs). Here, we introduce a sequential data assimilation method for incorporating observational data in a shape-morphing solution (SMS). Our method takes the form of a predictor-corrector scheme, where the observations are used to correct the SMS parameters using Newton-like iterations. Between observation points, the SMS equations---a set of ordinary differential equations--- are used to evolve the solution forward in time. We prove that, under certain conditions, the data assimilated SMS (DA-SMS) converges uniformly towards the true state of the system. We demonstrate the efficacy of DA-SMS on three examples: the nonlinear Schr\"odinger equation, the Kuramoto--Sivashinsky equation, and a two-dimensional advection-diffusion equation. Our numerical results suggest that DA-SMS converges with relatively sparse observations and a single iteration of the Newton-like method.
\end{abstract}

%
%

\section{Introduction}\label{sec:intro}
Shape-morphing solutions (SMS) are a new class of approximate solutions to spatiotemporal partial differential equations (PDEs). A distinguishing feature of SMS is that they depend nonlinearly on a set of time-dependent parameters. In their simplest form, SMS can be expressed as
\begin{equation}\label{eq:sms}
\hat{u}(\vc x, \params(t)) = \sum_{i=1}^\nummodes a_i(t) \phi_i(\vc x, \vc \beta_i(t)),
\end{equation}
where the parameters $\params(t) = \{a_i(t),\vc\beta_i(t)\}_{i=1}^\nummodes$ comprise the amplitudes $a_i$ and the shape parameters $\vc\beta_i$.
Equation~\eqref{eq:sms} can be viewed as a spectral method with the basis functions $\phi_i$~\cite{Anderson2022a,Anderson2024a}. Unlike classical spectral methods where the basis functions are static, shape-morphing solutions allow the basis to change over time via their dependence on the shape parameters $\vc\beta_i(t)$. As such, shape-morphing solutions can more easily adapt to the true solution of the PDE.

An alternative interpretation is to view equation~\eqref{eq:sms} as a shallow neural network with time-dependent parameters. In this case, $\phi_i$ is the activation function of the $i$-th node, $\vc\beta_i$ are the weights and biases of the node, and $a_i$ is the corresponding output weight. This interpretation can be easily extended to deep neural networks with time-dependent parameters~\cite{Du2021}.

Given a shape-morphing solution $\hat u(\vc x,\params(t))$, the evolution of its parameters $\params(t)$ are determined by a set of ordinary differential equations (ODEs), ensuring that the SMS closely approximates the true solution of the PDE. As we review in Section~\ref{sec:relWork}, the theory of shape-morphing solutions has been developing rapidly over the past few years. In the literature, the approximate solution $\hat u(\vc x,\params(t))$ is referred to by various names, such as evolutional deep neural networks~\cite{Du2021},
reduced-order nonlinear solutions~\cite{Anderson2022a}, and neural Galerkin schemes~\cite{Bruna2024}. Here, with hopes of establishing a unifying terminology, we refer to them as shape-morphing solutions or SMS, for short. We refer to the ODEs that govern the evolution of the parameters $\params(t)$ as the shape-morphing equation or the SMS equation.

The purpose of the present work is to introduce a data assimilation method which uses observational data from the physical system to steer the shape-morphing solution closer to the true dynamics. Data assimilation is generally used to alleviate various sources of error, such as modeling error, truncation error, and uncertainties in the initial condition, which inevitably affect any numerical forecast~\cite{Law2015}. Our method takes the form of a predictor-corrector scheme, where the shape-morphing equations are used to evolve the solution between observation points. At time instances where observational data are available, they are utilized to correct for errors using Newton-like iterations. We prove that, if enough observations are available, the data assimilated SMS converges uniformly towards the ground truth. We demonstrate the efficacy of our method on three numerical examples.

\subsection{Related Work}\label{sec:relWork}
Shape-morphing solutions were first developed by Du and Zaki~\cite{Du2021} in the context of deep neural networks and simultaneously by Anderson and Farazmand~\cite{Anderson2022a} in the context of model order reduction. Du and Zaki~\cite{Du2021} approximate the solution of a PDE using a deep neural network whose weights and biases are time-dependent, hence referring to the resulting method as the \emph{evolutional deep neural network} (EDNN). They derive a set of ordinary differential equations for evolving the parameters of the network, eliminating the need for a training process.
Anderson and Farazmand~\cite{Anderson2022a} independently arrived at the same set of equations where they considered reduced-order models whose modes depend nonlinearly on a set of time-dependent parameters, thus allowing for the modes to adapt to the solution of the PDE. The resulting method is referred to as \emph{reduced-order nonlinear solutions}  (RONS).

Ever since these original papers, there has been a rapidly growing body of work addressing various theoretical and computational aspects of shape-morphing solutions.
For instance, Bruna et al.~\cite{Bruna2024} introduced an adaptive sampling strategy to reduce the computational cost of the Monte Carlo sampling required by EDNN. Coining yet another terminology, they refer to their method as \emph{neural Galerkin schemes}.

Anderson and Farazmand~\cite{Anderson2024a} develop a number of methods to address the computational cost of SMS: (i) They developed a symbolic computing approach to reduce the computational cost of forming the SMS equations. (ii) When symbolic computing is not feasible, they introduced a collocation point method as an alternative and uncovered its relation to the Monte Carlo sampling method of~\cite{Du2021}. (iii) They also introduced a regularization to the SMS equations which speeds up the time integration. We also refer to Ref.~\cite{Finzi2023} which addresses similar computational issues. Anderson and Farazmand~\cite{Anderson2024b} applied these computational advances to solve the Fokker--Planck equation in higher dimensions and discovered that the metric tensor in SMS coincides with the Fisher information metric.

Shape-morphing solutions may over-parameterize the solution of the PDE, leading to redundant parameters. Berman and Peherstorfer~\cite{Berman2023} introduced a randomized projection method for eliminating such redundant parameters.

Enforcing the PDE's boundary conditions in shape-morphing solutions, and in particular evolutional neural networks, can be challenging~\cite{Du2021,Sukumar2022}. Kast and Hesthaven~\cite{Hesthaven2024} use embedding methods to ensure that the boundary conditions of the PDE are enforced directly into the neural network architecture.

Anderson and Farazmand~\cite{Anderson2022a} had already derived a set of SMS equations which enforce nonlinear constraints such as conservation of mass, energy, and other first integrals.
As mentioned by Hilliard and Farazmand~\cite{Hilliard2024}, these equations conserve the first integrals in continuous-time, necessitating a special time stepping scheme to ensure their conservation after discretization. To address this issue, Schwerdtner et al.~\cite{Schwerdtner2024} propose a discrete-in-time version of conservative SMS.

Shape-morphing solutions can be viewed as the flow of a vector field on a manifold embedded in the function space of the PDE~\cite{Anderson2022a}. Often this vector field is evaluated pointwise during the time integration process. Gaby et al.~\cite{Gaby2024} propose to learn the vector field on the SMS manifold as a neural network. They subsequently evolve the SMS equations by replacing the vector field with its neural net representation.

Following Ref.~\cite{Du2021}, Kim and Zaki~\cite{Zaki2024} developed a multi-network approach to shape-morphing solutions so that each state variable is represented by a separate EDNN. In addition, they partition the spatial domain into subsets so that the PDE is solved on each subset with its own corresponding EDNN.

In spite of the fast growing body of work on this topic, assimilating observational data in shape-morphing solutions has remained unaddressed. The purpose of the present paper is to introduce a data assimilation (DA) method for SMS approximation of PDEs. To this end, standard data assimilation methods, such as extended Kalman filters~\cite{Kalman1960, Kalman1961,McElhoe1966,Schmidt1966}, variational DA~\cite{Andersson1998,Courtier1998,Rabier1998}, Bayesian DA~\cite{Apte2007,reich2015}, and nudging methods~\cite{Ledimet1992,Azouani2014} can in principle be used. However, compared to our proposed method, Kalman filtering and variational methods are computationally more expensive as they require multiple solves of the underlying PDE and/or the corresponding adjoint equation~\cite{Tremolet2006}. Furthermore, the nonlinear dependence of SMS $\hat u(\vc x,\params(t))$ on its parameters $\params$ inhibits straightforward implementation of classical DA methods. In contrast, our proposed predictor-corrector DA scheme is easy to implement and avoids repeated solves of the PDE.

\subsection{Main Contributions}
The main contributions of this paper can be summarized as follows.
\begin{enumerate}
	\item Discrete-time data assimilated SMS (DA-SMS): We introduce a sequential data assimilation algorithm that propagates the solution forwards in time and incorporates observational data into the solution as they becomes available at discrete times.
	\item Convergence result: We prove that, if enough observational data is available, the DA-SMS method converges towards the true solution of the system.
	\item Continuous-time data assimilation: We also present a DA method for SMS in case the observational data is measured with high temporal frequency.
	\item As a byproduct of our work, we introduce a new method for directly encoding mixed Dirichlet--Neumann boundary conditions in the architecture of a neural network.
	\item We demonstrate the application of our method on three examples with one and two spatial dimensions.
\end{enumerate}

\subsection{Outline}
The remainder of the paper is organized in the following manner. In Section~\ref{sec:prelim}, we review the method of shape-morphing solutions and how we evolve the solutions forward in time. In Section~\ref{sec:theory}, we introduce our main results for sequential data-assimilation. We discuss both discrete DA-SMS algorithm and its continuous-time version. In Section~\ref{sec:numerics}, we present three numerical examples: the nonlinear Sch\"{o}dinger (NLS) equation in Section~\ref{sec:NLS}, the Kuramoto--Sivashinsky (KS) equation in Section~\ref{sec:KS}, and the two-dimensional advection-diffusion (AD) equation in Section~\ref{sec:ADE}. We present our concluding remarks in Section~\ref{sec:conc}.

\section{Mathematical Preliminaries}\label{sec:prelim}
In this section, we summarize the theory of shape-morphing modes to obtain numerical solutions to PDEs~\cite{Anderson2022a, Anderson2024a,Anderson2024b,Anderson2022b}.  
Consider a PDE of the form,
\begin{equation}\label{eq:gen_PDE}
\pard{u}{t} = F(u),  \quad u(\vc x, 0) = u_0(\vc x),
\end{equation}
where $F$ is a differential operator and $u:\spatdom\times \R^+\rightarrow \R$ denotes the solution, where $\spatdom\subset \R^d$ is the spatial domain of the PDE. For any time $t$, the solution $u(\cdot,t)$ belongs to an appropriate Hilbert space $\hilspace$ with associated inner product $\hilip{\cdot}{\cdot}$ and induced norm $\hilnorm{\cdot}$.  For the purposes of this exposition, we assume that the solution $u(\vc x, t)$ is scalar-valued. However, this framework can easily be extended to vector-valued systems of PDEs as discussed in Ref.~\cite{Hilliard2024}. Additionally, here we only consider real-valued PDEs, but the extension to complex-valued PDEs is straightforward and can be found in Refs.~\cite{Anderson2022b,Hilliard2024}.  

Next, we introduce the approximate shape-morphing solution,
\begin{equation} \label{eq:gen_sms}
\hat{u}\left(\vc x,\params(t)\right)  = \sum_{i=1}^\nummodes a_i(t) \phi_i(\vc x, \bs \beta_i (t)),
\end{equation}
with the \textit{shape-morphing modes} $\phi_i$, the corresponding amplitudes $a_i$, and the vector of \textit{shape parameters} $\vc \beta_i\in\mathbb R^{\numparams-1}$. We denote the full set of parameters by $\params(t) = \{a_i(t),\bs \beta_i(t)\}_{i=1}^N\in \R^{pN}$, where each term in~\eqref{eq:gen_sms} contains $\numparams$ parameters.
The choice of the shape-morphing modes depends on the PDE. For instance, one can take a Gaussian mixture model where the shape-morphing modes are Gaussians and the shape parameters are the mean and the variance~\cite{Anderson2022b,Anderson2024b}. A more general choice is an artificial neural network (ANN), where $\phi_i$ is the activation function, and the shape parameters are the corresponding weights and biases of the network~\cite{Anderson2024a,Du2021}. For example, a shallow neural network with a hyperbolic tangent activation function coincides with the shape-morphing solution,
\begin{equation}
\hat{u}\left(x,\params(t)\right) = \sum_{i=1}^\nummodes a_i(t) \tanh\left(w_i(t) x + b_i(t)\right),
\end{equation}
whose shape-parameters are the $\nummodes$ amplitudes $a_i$, weights $w_i$, and biases $b_i$, so that $\params = \{a_i,w_i,b_i\}_{i=1}^N \in \R^{3N}$. 

In the conventional approach to ANNs for PDEs, the parameters are determined by minimizing a loss function~\cite{Karniadakis2021,Stuart2023,Karniadakis2021b}.
The crucial difference in our approach is that SMS eliminates the need for this training step by evolving the time-dependent parameters $\params(t)$ through an explicit set of ordinary differential equations (ODEs). 

To derive the ODEs governing the evolution of parameters, we introduce the residual function,
\begin{equation}\label{eq:residual}
R(\vc x, \params, \dot{\params}) := \hat{u}_t - F(\hat{u}) = \sum_{j=1}^{\numparams N} \pard{\hat{u}}{\singparam_j}\dot{\singparam}_j - F(\hat{u}). 
\end{equation}
which measures the discrepancy between the rate of change of the approximate shape-morphing solution and the dynamics governed by the PDE. To derive the corresponding system of ODEs, we seek $\dot{\params}$ which minimizes the instantaneous error, 
\begin{equation}\label{eq:sym_opt}
\min_{\dot{\params}\in \R^{\numparams N}} \frac12\hilnorm{R(\cdot, \params,\dot{\params})}^2. 
\end{equation}
The solution to this minimization problem is given explicitly in terms of the \textit{shape-morphing equation}~\cite{Anderson2022a,Du2021},
\begin{equation}\label{eq:SME_sym}
M(\params)\dot{\params} = \vc f(\params),
\end{equation}
where $M\in \R^{\numparams N\times \numparams N}$ is the \textit{metric tensor},
\begin{equation}
M_{ij}(\params) = \hilip{\pard{\hat{u}}{\singparam_i}}{\pard{\hat{u}}{\singparam_j}},
\end{equation}
and the components of the vector field $\vc f:\mathbb R^{pN}\to\mathbb R^{pN}$ are given by
\begin{equation}
f_i(\params) = \hilip{\pard{\hat{u}}{\singparam_i}}{F(\hat{u})}. 
\end{equation}
This set up is especially useful in the event that one can symbolically compute the integrals in Hilbert inner products for $M$ and $\vc f$. In this case, one can exploit the structure of the metric tensor and the vector $\vc f$ to make the implementation of symbolic computations more efficient. We refer to Ref.~\cite{Anderson2024a} for further detail regarding these computational aspects.

When symbolic computations of $M$ and $\vc f$ are not feasible, we instead opt for the collocation point method which was originally derived in Ref.~\cite{Anderson2024a}. In this case, we specify the $\numcoll$ collocation points $\{\vc x_1, \ldots, \vc x_{\numcoll}\}\subset\Omega$, and require that the residual function vanishes at these points, i.e. 
\begin{equation}\label{eq:1_coll}
R(\vc x_i, \params, \dot{\params}) = \sum_{j=1}^{\numparams N}  \pard{\hat{u}}{\singparam_j}(\vc x_i,\params)\dot{\singparam}_j - F(\hat{u})\bigg\rvert_{\vc x = \vc x_i}=0, \quad  i = 1, \ldots, \numcoll.
\end{equation}
We can more conveniently express \eqref{eq:1_coll} in one equation by defining the \textit{collocation metric tensor} $\tilde{M}(\params)\in \R^{\numcoll\times \numparams N}$ whose elements are give by
\begin{equation}
\tilde{M}_{ij} (\params)= \pard{\hat{u}}{\singparam_j}(\vc x_i, \params) ,
\end{equation}
and the vector field $\tilde{\vc f}\in \R^{\numcoll}$ such that
\begin{equation}
\tilde{f}_i(\params) = F\left(\hat{u}(\cdot, \params)\right)\bigg\rvert_{\vc x = \vc x_i}.  
\end{equation}
With these definitions, we can express \eqref{eq:1_coll} as
\begin{equation}\label{eq:coll_sms}
\tilde{M}(\params) \dot{\params} = \tilde{\vc f}(\params).
\end{equation}
Note that no integration is required in the construction of the collocation metric tensor $\tilde{M}$ or the vector field $\tilde{\vc f}$.

The collocation metric tensor $\tilde{M}$ is rectangular which means that the system \eqref{eq:coll_sms} can be either overdetermined or underdetermined. In such case, we obtain $\dot{\params}$ by solving the linear least squares problem,  
\begin{equation} \label{eq:coll_opt}
\min_{\dot{\params}\in \R^{\numparams N}}\frac{1}{2}\|\tilde{M}(\params)\dot{\params} - \tilde{\vc f}(\params)\|_2^2, 
\end{equation}
which is equivalent to minimizing $\sum_i |R(\vc x_i,\params,\dot{\params})|^2$. Here, $\|\cdot\|_2$ denotes the Euclidean norm.
The resulting \textit{collocation shape-morphing equation} is
\begin{equation}\label{eq:coll_SMS}
\dot{\params} = \tilde{M}^+(\params) \tilde{\vc f}(\params),
\end{equation}
in terms of the Moore--Penrose pseudo-inverse $\tilde{M}^+$.

Lastly, in either the shape-morphing equation \eqref{eq:SME_sym} or the collocation shape-morphing equation \eqref{eq:coll_sms}, the metric tensors can be ill-conditioned which leads to stiffness in the corresponding ODEs. In order to alleviate this issue, Anderson and Farazmand~\cite{Anderson2024a} regularize the optimization problems~\eqref{eq:sym_opt} and \eqref{eq:coll_opt} using Tikhonov regularization. Introducing the regularization parameter $\SMSTikh$, the regularized optimization problem~\eqref{eq:sym_opt} reads
\begin{equation}
\min_{\dot{\params} \in \R^{\numparams N}} \frac{1}{2}\hilnorm{R(\vc x, \params, \dot{\params})}^2 + \frac{\SMSTikh}{2} \|\dot{\params}\|_2^2,
\end{equation}
whose solution is given by
\begin{equation}\label{eq:reg_sym}
\left(M(\params)+ \SMSTikh I \right) \dot{\params} =  \vc f(\params),
\end{equation}
where $I$ denotes the $\numparams N\times \numparams N$ identity matrix. We similarly regularize the least squares problem~\eqref{eq:coll_opt} and solve
\begin{equation}
\min_{\dot{\params}\in \R^{\numparams N}}\frac{1}{2}\|\tilde{M}(\params)\dot{\params} - \tilde{\vc f}(\params)\|_2^2 + \frac{\SMSTikh}{2} \|\dot{\params}\|^2_2.
\end{equation}
The resulting regularized collocation shape-morphing equation reads
\begin{equation}\label{eq:reg_coll}
\left(\tilde M^\top(\params) \tilde M(\params) + \SMSTikh I \right) \dot{\params} = \tilde M^\top(\params) \tilde{\vc f}(\params).
\end{equation}
To evolve the parameters $\params$ in time, we solve the ODEs \eqref{eq:reg_sym} or \eqref{eq:reg_coll} using standard numerical time integrators such as explicit Runge--Kutta schemes.

Lastly, we emphasize that the metric tensors $M$ and $\tilde{M}$ as well as the vectors $\vc f$ and $\tilde{\vc f}$ are functions of the parameters $\params$; but for the remainder of this manuscript we omit their dependence on $\params$ for notational simplicity. 

\section{Data Assimilation}\label{sec:theory}
In this section, we introduce our sequential data assimilated shape-morphing solutions algorithms. More specifically, in Section \ref{sec:DA-SMS} we present the main result of our paper: a discrete sequential data assimilation algorithm used in conjunction with the shape-morphing solutions. We refer to this methodology as data-assimilated SMS or DA-SMS for short. In Section \ref{sec:sens_place}, we derive a convergence guarantee for DA-SMS. Lastly, in Section \ref{sec:cont_DA}, we introduce a continuous time version of DA-SMS that instead uses the time derivative of the observational data. 
 
\subsection{Discrete Sequential Data Assimilation: DA-SMS}\label{sec:DA-SMS}
Consider a system described by PDE \eqref{eq:gen_PDE} with the initial condition $u_0$. Any method for approximating the system's evolution
needs to contend with various sources of error: modeling error resulting from simplifying assumptions in deriving the PDE, uncertainties in the initial condition, truncation error from representing the solution with a finite-dimensional approximation, and time integration error. Data assimilation seeks to remedy these errors by incorporating direct observational data from the system into the numerical solution. Here, we develop data assimilation methods that are applicable when the approximate solution is obtained as a shape-morphing solution $\hat{u}(\cdot,\params(t))$ with nonlinear dependence on the parameters $\params$.

To begin, we assume that we have access to $\numsens$ sensors which take measurements of observable quantities associated with the system. We denote these observables by $y_j(t)$ for $j\in \{1,\ldots, \numsens\}$. These measurements are collected at $m$ discrete points in time: $t_1< t_2< \ldots< t_m$. We refer to the time interval between the initial time $t_0$ and the final observation time $t_m$ as the \textit{data-assimilation (DA) time window}, $\mathcal{T}_{DA} = [t_0,t_{\numtimes}]$.  At each $t_i$, we collect the $\numsens$ measurements into a vector $\vc \data(t_i) \in \R^\numsens$ whose components we denote by, 
\begin{equation}
\data_j(t_i) = C_j(u(\cdot, t_i)) + \eta, \quad j = 1,\ldots, \numsens,
\end{equation}
where $C_j:\mathcal{H} \rightarrow \R$ is the associated observation operator and $\eta$ is the observational error or noise. The observables can be pointwise observations of the state of the system, i.e. $y_j(t_i) = u(\vc x_j, t_i)$, or they could be a nonlinear operator applied to the state of the system such as the modulus of the solution to the nonlinear Schr\"{o}dinger equation as discussed later in Section \ref{sec:NLS}. 
\begin{figure}[t]
\centering 
\includegraphics[width = .85\textwidth]{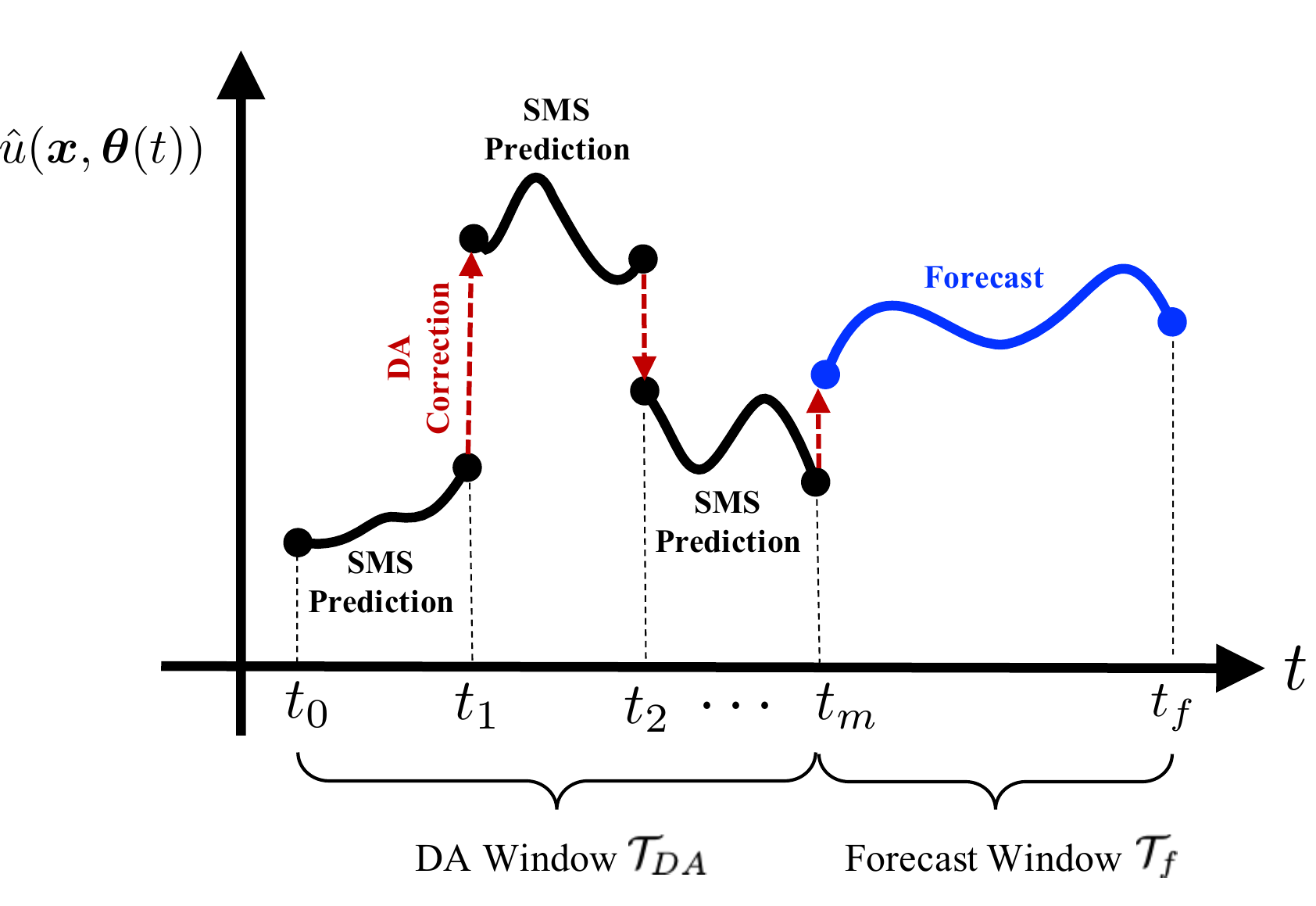}
\caption{\label{fig:DA_diagram}Schematic illustration of the DA-SMS algorithm.}
\end{figure}

Figure \ref{fig:DA_diagram} illustrates the overall DA-SMS algorithm. Our goal is to use the data to correct for the errors introduced by the approximations. Our algorithm achieves this objective in two main steps. The first is the prediction step where we use the SMS equations derived in Section \ref{sec:prelim} to evolve the approximate solution from time $t_{i-1}$ to $t_{i}$. In Figure \ref{fig:DA_diagram}, this step is represented by the black curves. The second piece of the DA-SMS algorithm, is the correction step, depicted by the red dashed lines in Figure \ref{fig:DA_diagram}. At time $t_i$, we freeze time and make a correction to the parameters $\params(t_i)$, using the true observational data $\vc y(t_i)$. 

More specifically, during the correction step, we compute the SMS estimated observations defined by,  
\begin{equation}
\hat{\data}_j(t_i) = C_j[\hat{u}(\cdot, \params(t_i))], \quad j = 1,\ldots, \numsens.
\end{equation}
For notational simplicity, we often write $C_j(\params(t_i))$ instead of $C_j[\hat{u}(\cdot, \params(t_i))]$. In the DA step of our algorithm, we seek a correction to the parameters $\params(t_i)$ that minimizes the discrepancy between the true data $\vc \data(t_i)$ of the system and the predicted data $\hat{\vc \data}(t_i)$. To this end, we solve the nonlinear equation,
\begin{equation}\label{eq:NLSQ}
\vc C(\params(t_i)+\delta\params) = \vc y(t_i),
\end{equation} 
where $\delta\params$ is the correction step and $\vc C = [C_1(\cdot), C_2 (\cdot), \ldots, C_r(\cdot)]^\top$. Note that even if the observation operator is linear, $\vc C(\params(t_i))$ is a nonlinear function of $\params(t_i)$ due to the nonlinear dependence of the shape-morphing solution $\hat{u}(\cdot, \params)$ on its parameters.

To solve the nonlinear equation~\eqref{eq:NLSQ}, we introduce the Newton-like iterations,
\begin{equation}\label{eq:newton-like}
\begin{split}
\params^{(k+1)}(t_i) & = \params^{(k)}(t_i)+\left[\nabla_{\params}\vc C\left(\params^{(k)}(t_i)\right)\right]^+ \left[\boldsymbol \data(t_i)-\vc C\left(\params^{(k)}(t_i)\right) \right], \\
\end{split} 
\end{equation}
where the superscript on $\params$ denotes the $k$-th iteration. 
\begin{table}
\centering 
\begin{tabular}{| p{.15 \textwidth} | p{.75\textwidth} |}
\hline 
\textbf{Symbol} & \textbf{Description} \\
\hline
$\params_i^{(k)}$ & Parameters at time $t_i$ and the $k$-th Newton iteration: $\params_i^{(k)} = \params^{(k)}(t_i)$ \\ 
\hline
$\hat{\vc y}_i^{(k)}$ & SMS estimated observations at time $t_i$ and $k$-th iteration: $\hat{\vc y}_i^{(k)} = \vc C[\hat{u}(\cdot, \params_i^{(k)})]$\\ 
\hline
$\vc y_i$ & True observational data at time $t_i$: $\vc y_i = \vc y(t_i)$ \\
\hline
\end{tabular}
\caption{\label{tbl:notation} The quantities used in the DA-SMS algorithm for notational simplicity.}
\end{table}

To further simplify the notation, we introduce the following. We denote the parameters at time $t_i$ and the $k$-th Newton iteration by $\params_i^{(k)}=\params^{(k)}(t_i)$, the true observational data at time $t_i$ by $\vc \data_i=\vc \data(t_i)$, and the SMS estimated data at time $t_i$ and the $k$-th Newton iteration by $\hat{\vc y}^{(k)}_i=\vc C[\hat{u}(\cdot, \params^{(k)}(t_i))]$. In each of these, the subscript denotes the $i$-th time instance, and the superscript denotes the $k$-th Newton iteration. We summarize this notation in Table \ref{tbl:notation}. Lastly, we introduce the Jacobian matrix $J_r(\params) = \nabla_{\params}C(\params)\in \R^{r\times \numparams N}$, so that
\begin{equation}\label{eq:jac_def}
\jacop\left(\params\right) = \begin{bmatrix}
\mbox{---}\ \left[\nabla_{\params} C_1\left(\params\right)\right]^\top\ \mbox{---} \\
\mbox{---}\ \left[\nabla_{\params} C_2\left(\params\right)\right]^\top\ \mbox{---} \\
\vdots \\
\mbox{---}\ \left[\nabla_{\params} C_r\left(\params\right)\right]^\top\ \mbox{---}
\end{bmatrix}.
\end{equation} 

Using this notation, the iterates in \eqref{eq:newton-like} read
\begin{equation}\label{eq:DA_iterates}
\params^{(k+1)}_i = \params_i^{(k)} + \left[\jacop\left(\params_i^{(k)}\right)\right]^+ \left( \boldsymbol \data_i - \hat{\vc \data}_i^{(k)} \right).
\end{equation}
We stop the iterations once the relative error between the true observations, $\vc \data_i$, and the SMS estimated observations $\hat{\vc y}^{(k)}_i$ falls below a user-defined tolerance $\delta$, 
\begin{equation}
\frac{\|\hat{\vc \data}^{(k)}_i - \vc \data_i\|_2}{\|\vc \data_i\|_2}<\delta . 
\end{equation}
The tolerance $\delta$ should be no less than the amount of observational noise in the data to avoid overfitting to noisy data. In addition, we also specify a maximum number of iterations at each data assimilation step to avoid infinite iterations. However, as we will see in Section \ref{sec:numerics}, often only one Newton iteration is sufficient in practice. 

In practice, the number of sensors or data points $\numsens$ is much smaller than the total number of parameters $\numparams \nummodes$. As such, the nonlinear problem in \eqref{eq:NLSQ} is underdetermined, and therefore its solution is not necessarily unique. Nonetheless, under certain conditions, underdetermined Newton iterations are guaranteed to converge to a solution \cite{Chen1994,Kelley1999}. In particular, underdetermined Newton-like iterations~\eqref{eq:DA_iterates} converge to a solution $\params_i^*$ under the following assumptions: 
\begin{enumerate}
\item The observation operator $\vc C(\params)$ is Lipschitz continuous in a neighborhood around $\params_i^*$, and
\item The Jacobian matrix $\jacop(\params_i^*)$ has full row rank, i.e., $\mbox{rank}(\jacop(\params_i^*)) = r$,
\end{enumerate}
as long as the initial guess $\params_i^{(0)}$ is close enough to a solution $\params_i^*$; see Theorem 2.4.2 of Ref.~\cite{Kelley1999}.

If the Jacobian matrix $\jacop$ is ill-conditioned, the Newton-like iterations \eqref{eq:DA_iterates} will amplify the observational noise. To address this issue, we use a regularized version of the Newton-like iterations~\cite{Pes2020,Rieder1999}. Letting $\DATikh$ be the regularization parameter, the iterates in \eqref{eq:DA_iterates} become, 
\begin{equation}\label{eq:reg_newton}
\begin{split}
\params_i^{(k+1)} & = \params_i^{(k)}+\jacop\left(\params_i^{(k)}\right)^\top  \left(\jacop\left(\params_i^{(k)}\right)\jacop\left(\params_i^{(k)}\right)^\top + \DATikh I\right)^{-1} \left( \boldsymbol \data_i-\hat{\vc \data}^{(k)}_i\right) \\
& = \params_i^{(k)}+\tilde{\jacop}\left(\params_i^{(k)}\right)^+ \left( \boldsymbol \data_i - \hat{\vc \data}_i^{(k)}\right),
\end{split}
\end{equation}
where $\tilde{\jacop}\left(\params_i^{(k)}\right)^+ = \jacop\left(\params_i^{(k)}\right)^\top  \left(\jacop\left(\params_i^{(k)}\right)\jacop\left(\params_i^{(k)}\right)^\top + \DATikh I\right)^{-1}$ is the regularized pseudo-inverse of $\jacop\left(\params_i^{(k)}\right)$.
%
%
\begin{algorithm}[t!]
	\caption{Data assimilated shape-morphing solutions (DA-SMS) algorithm.}\label{alg:DA-SMS} 
	\KwInput{Initial parameter values $\params_0$, Observation times $t_1<t_2< \ldots< t_{m}$, Observations $\{\vc \data_1, \ldots, \vc \data_{m}\}$, Maximum Newton iterations (maxits), Newton iteration tolerance $\delta$, Final forecast time $t_f$}
\For{$i = 1,2,\ldots, m$}
{
   $\params_{i}^{(0)}= \texttt{evolve}(\params_{i-1}, [t_{i-1},t_{i}])$ \Comment*[r]{Use SMS to predict from time $t_{i-1}$ to $t_i$}
   $\hat{\vc \data}_{i}^{(0)} = \vc C(\params_{i}^{(0)})$\;	
   $error=\| \vc \data_{i} - \hat{\vc \data}^{(0)}_{i} \|_2/ \|\vc \data_i \|_2$\;
   $k=0$\;
	\While {($error\geq \delta$ and $k\leq\text{maxits}$)}
	{
        $J_r =$ $\nabla_{\params} \vc C (\params_{i+1}^{(k)})$\;   
		$\params_i^{(k+1)} =\params_i^{(k)}+ \tilde{J_r}\left(\params_i^{(k)}\right)^+ (\vc \data_{i} - \hat{\vc \data}_{i}^{(k)})$ \Comment*[r]{Regularized Newton iteration}
		$\hat{\vc \data}^{(k+1)}_i = \vc C(\params_i^{(k+1)})$\; 	
		$error = \| \vc \data_{i} - \hat{\vc \data}^{(k+1)}_i \|_2/ \|\vc \data_i \|_2$\;
		$k = k+1$\;
	}
	$\params_{i} = \params_{i}^{(k)}$\;
}
$\params_f = $ \texttt{evolve}$(\params_{m},[t_m,t_f])$   \Comment*[r]{Use SMS to forecast to final time $t_f$}
\KwOutput{DA-SMS forecasted parameters $\params_{f}$}
\end{algorithm}

At time $t_m$, the last time where observational data is available, we obtain a final set of corrected parameters $\params_m$. Since no more data is available, we use the shape-morphing equations to evolve from $t_m$ to the final time $t_f$. We refer to this time interval as the \textit{forecast window} $\mathcal{T}_f= [t_m,t_f]$. At time $t_f$, we use the set of parameters $\params_f$ to reconstruct the final forecast of the system $\hat{u}(\cdot, \params_f)$. We summarize the DA-SMS method in Algorithm~\ref{alg:DA-SMS}.

\subsection{Convergence Analysis}\label{sec:sens_place}
Consider pointwise observations where $C_j[u(\cdot,t_i)]=u(\vc x_j,t_i)$ with $\vc x_j$ being a sensor location. For simplicity, we assume that the observations are noise-free; the results corresponding to noisy data are similar, as discussed in Remark~\ref{rem:noise_ub} below.
Convergence of the Newton iterations in Algorithm~\ref{alg:DA-SMS} ensures that the shape-morphing solution $\hat u(\cdot,\params_i)$ converges to the true solution $u(\cdot,t_i)$ at the sensor locations. However, this does not immediately guarantee convergence to the true solution at other points $\vc x \in\Omega$ where no sensors are available. In this section, we prove that, under certain conditions, the shape-morphing solution in fact converges uniformly to the true solution throughout the spatial domain $\Omega$.

Since time is frozen during the Newton iterations, we omit the dependence of the true solution $u(\vc x, t_i)$ on time and simply write $u(\vc x)$ throughout this section. In addition, we assume that $\spatdom\subset \R^n$ is bounded. Consider $\numsens$ sensors which are located at $\{\vc x_1, \ldots, \vc x_\numsens\}\subset\Omega$, so that $u(\vc x_j)$ is known for all $j \in \{1,2,\ldots, r\}$. We define
\begin{equation}
\Delta = \sup_{\vc x\in \spatdom} \inf_{1\leq j \leq \numsens}\|\vc x - \vc x_j \|,
\end{equation}
which quantifies the maximum distance between any point in $\Omega$ and its closest sensor. More specifically, for any $\vc x\in \spatdom$, the ball $B_{\Delta}(\vc x)$, with radius $\Delta$ centered at $\vc x$, contains at least one sensor. 
\begin{ass}\label{ass:sens_place} We assume that the true solution $u$ and the shape-morphing approximation $\hat{u}$ have the following properties.
\begin{enumerate}
\item Both $u$ and $\hat{u}$ are Lipschitz continuous with respect to $\vc x$. We denote their smallest Lipschitz constants by $L_{u}$ and $L_{\hat{u}}$ respectively. 
\item There exists a sequence of parameters $\{\params^{(1)}, \params^{(2)}, \ldots\}$ such that 
\begin{equation}\label{eq:pw_conv}
\lim_{k\rightarrow\infty} \big | \hat{u}(\vc x_j, \params^{(k)}) - u(\vc x_j) \big| = 0, \quad \forall j\in\{1,\ldots, r\}.
\end{equation}
\end{enumerate}
\end{ass}
These assumptions are not necessarily strong constraints on our approximate solution $\hat{u}$. Many choices of shape morphing solutions, such as Gaussian mixture models or neural networks, are Lipschitz continuous~\cite{GoodBengCour16,Virmaux2018}. Additionally, the second assumption simply ensures that our approximate solution $\hat{u}$ is capable of approximating the unknown function $u$ at the sensor locations. This assumption can be easily satisfied by choosing shape-morphing modes which satisfy the universal approximation theorem. In equation~\eqref{eq:pw_conv}, we used an arbitrary convergent sequence of parameters $\params^{(k)}$. Although in this paper we obtain this sequence through Newton-like iterations~\eqref{eq:DA_iterates}, alternative methods can also be used. 

The following lemma provides an upper bound for the pointwise error between the true solution $u$ and the shape-morphing approximation $\hat u$.
\begin{lem} \label{lem:sens_place}
Let Assumption~\ref{ass:sens_place} hold. For any $\vc x\in \spatdom$, let $\hat{j}$ denote the integer such that $\vc x_{\hat{j}}$ is the closest sensor to point $\vc x$. Then for any parameter $\params$, we have 
\begin{equation}
|\hat{u}(\vc x, \params) - u(\vc x)| \leq \left(L_{u} + L_{\hat{u}}\right)\Delta + |\hat{u}(\vc x_{\hat{j}},\params) - u(\vc x_{\hat j}))|. 
\end{equation}
\begin{proof}
Using the triangle inequality, we have 
\begin{align*}
|\hat{u}(\vc x, \params) - u(\vc x)| & = |\hat{u}(\vc x, \params)- \hat{u}(\vc x_{\hat{j}},\params) + \hat{u}(\vc x_{\hat{j}},\params) - u(\vc x_{\hat j}) + u(\vc x_{\hat j}) - u(\vc x)|, \\
& \leq |\hat{u}(\vc x, \params)- \hat{u}(\vc x_{\hat{j}},\params)|  + |\hat{u}(\vc x_{\hat{j}},\params) - u(\vc x_{\hat j})| + |u(\vc x_{\hat j}) - u(\vc x)|,\\
& \leq (L_{u}+ L_{\hat{u}})\|\vc x - \vc x _{\hat{j}}\| + |\hat{u}(\vc x_{\hat{j}}, \params) - u(\vc x_{\hat{j}})|,
\end{align*}
where we used the Lipschitz continuity of the functions for the last inequality.
Since $\vc x_{\hat j}$ is the closest sensor to $\vc x$, we have that $\|\vc x - \vc x_{\hat{j}}\| \leq \Delta$, which concludes the proof. 
\end{proof}
\end{lem}

We now state our main result which shows that, if enough sensors are provided, the shape-morphing approximation $\hat u$ becomes arbitrarily close to the true solution $u$ in the uniform norm.
\begin{thm} \label{thm:sens_place}
Let Assumption~\ref{ass:sens_place} hold. Furthermore, for any $\epsilon > 0$, assume that there are enough sensors distributed such that $\Delta < \epsilon/2(L_{u} + L_{\hat{u}})$. Then there exists $K\in \N$ such that for any $k> K$, 
\begin{equation}\label{eq:ub_sens_place}
\sup_{\vc x\in \spatdom} |\hat{u}(\vc x, \params^{(k)}) - u(\vc x)|<\epsilon.
\end{equation}
\end{thm}
\begin{proof}
	Consider the upper bound in Lemma~\ref{lem:sens_place}. Equation~\eqref{eq:pw_conv} in Assumption~\ref{ass:sens_place} implies that $\lim_{k\rightarrow \infty} |\hat{u}(\vc x_{\hat{j}},\params^{(k)})-u(\vc x_{\hat{j}})| = 0$. Therefore, there exists $K\in \N$, such that for any $k>K$ we have $|\hat{u}(\vc x_{\hat{j}},\params^{(k)})-u(\vc x_{\hat{j}})|<\epsilon/2$. Moreover, since Lemma \ref{lem:sens_place} holds for any $\params$, it also holds for the sequence $\{\params^{(k)}\}_{k\geq1}$. Using this fact, Lemma \ref{lem:sens_place}, and the fact that $\Delta < \epsilon/2(L_{u}+L_{\hat{u}})$, we have
\begin{equation}
|\hat{u}(\vc x, \params^{(k)})-u(\vc x)|<\epsilon,
\end{equation}
for any $\vc x\in \spatdom$ and $k>K$. 
\end{proof}

This result yields both an indication of the number of required sensors as well as their distribution throughout the domain $\Omega$. However, it is important to note that these bounds are rather pessimistic. As we show in Section~\ref{sec:numerics}, a relatively smaller number of sensors, uniformly distributed over the spatial domain, is often sufficient in practice. 

\begin{rem}\label{rem:noise_ub}
In proving Theorem~\ref{thm:sens_place}, we assumed that the observational data is noise-free. In the realistic case, where observational data $u(\vc x_j)+\eta$ is polluted with the noise $\eta$, equation~\eqref{eq:pw_conv} needs to be modified to read, $\lim_{k\to\infty} |\hat u(\vc x_j,\params^{(k)})-u(\vc x_j)|\leq |\eta|$. Subsequently, assuming that 
$\eta\sim \mathcal N(0,\sigma^2)$, the upper bound in~\eqref{eq:ub_sens_place} changes to
\begin{equation}
	\sup_{\vc x\in \spatdom}\mathbb E\left[ |\hat{u}(\vc x, \params^{(k)}) - u(\vc x)|\right]<\epsilon+\sqrt{\frac{2}{\pi}}\sigma,
\end{equation}
where $\mathbb E[\cdot]$ denotes the expected value. Although $\epsilon$ can be made arbitrarily small by increasing the number of sensors and the iteration number $k$, the error upper bound remains positive due to the variance $\sigma^2$ of the observational noise.
\end{rem}

\subsection{Continuous-time Data Assimilation}\label{sec:cont_DA}
The algorithm developed in Section~\ref{sec:DA-SMS} considers data obtained at discrete time instances, which is the situation one encounters in practice. Here, for completeness, we also introduce a data assimilation method for shape-morphing modes in case observational data are available in continuous time. To begin, we consider the instantaneous optimization problem \eqref{eq:sym_opt} from which we derived the shape-morphing equation~\eqref{eq:SME_sym}. The extension of our framework to the collocation shape-morphing equation~\eqref{eq:coll_SMS} is straightforward, as shown in~\ref{app:cont_deriv}. 

Recall from Section \ref{sec:prelim} that we seek to minimize the residual function $R(\cdot,\params, \dot{\params})$ which quantifies the discrepancy between the rate of change $\hat{u}_t$ of the approximate solution and the right-hand side of the PDE $F(\hat{u})$. In this case we assume that we have access to a continuum of data $\vc \data(t)$, and that we can compute its time derivative $\dot{\vc \data}$. To assimilate this data, we require that the time derivative of the SMS estimated observations $\dot{\hat{\vc \data}}$ are equal to $\dot{\vc \data}$, i.e. $\dot{\hat{\vc \data}} = \dot{\vc \data}$. In terms of the observation operators $C_j$, we require, 
\begin{equation}\label{eq:obs_dt}
\dot{\hat{ y}}_j = \frac{\id\;}{\id t} C_j( \params)  = \left\langle\nabla_{\params}[C_j (\params)], \dot{\params}\right\rangle = \dot{ \data}_j, \quad \forall j \in \{1,\ldots, \numsens\},
\end{equation}
where $\langle\cdot,\cdot\rangle$ denotes the Euclidean inner product and we used the fact that the observation operator $C_j$ has no explicit time dependence. 
Using the definition of $\jacop$ from \eqref{eq:jac_def}, we write equation~\eqref{eq:obs_dt} more compactly as
 \begin{equation}\label{eq:cont_DA_1}
\jacop(\params) \dot{\params} = \dot{\vc \data}. 
\end{equation}

To arrive at the \textit{continuous-time DA-SMS equation}, we constrain the optimization problem \eqref{eq:sym_opt} so that \eqref{eq:cont_DA_1} is satisfied. Namely, we seek the solution to
the constrained optimization problem,
\begin{equation}\begin{split}\label{eq:cont_DA_opt}
\min_{\dot{\params}} & \quad \frac{1}{2}\hilnorm{R(\cdot, \params, \dot{\params})}^2, \\ 
\text{s.t.} & \quad \jacop(\params) \dot{\params} = \dot{\vc \data}. 
\end{split}
\end{equation}
As we show in \ref{app:cont_deriv}, if the metric tensor $M$ is non-singular and $\jacop$ is full-rank, the solution to the constrained optimization problem \eqref{eq:cont_DA_opt} is given by
\begin{equation}\label{eq:cont_DA}
M\dot{\params} = \left( I -  \jacop^\top\left(\jacop  M^{-1} \jacop^\top \right)^{-1} \jacop M^{-1} \right) \vc f +  \jacop^\top \left(\jacop  M^{-1} \jacop^\top\right)^{-1}\dot{\vc y},
\end{equation}
where we have suppressed the dependence of $M$, $\vc f$, and $\jacop$ on $\params$. The solution to the system of ODEs in \eqref{eq:cont_DA} ensures that $\dot{\hat{\vc y}}(t) = \dot{\vc \data}(t)$ for all times $t$. 

We emphasize that, although this methodology can be useful theoretically, it is not necessarily suited for practical applications. First, all real-world observational data are polluted by noise; numerical computations of the time derivative will amplify this noise. In addition, if there is an initial difference in the SMS approximate observations, $\hat{\vc \data}(0) = \vc \data(0) + c$, then this difference will persist in the solution of~\eqref{eq:cont_DA}. Due to these issues, we do not recommend this methodology for most practical applications, and instead suggest using the discrete DA-SMS algorithm presented in Section~\ref{sec:DA-SMS}.

\section{Numerical Results}\label{sec:numerics}
We present three numerical examples: the nonlinear Schr\"{o}dinger (NLS) equation, the Kuramoto--Sivashinsky (KS) equation, and the advection-diffusion (AD) equation. For each example, we compare direct numerical simulations (DNS) to approximate solutions obtained from SMS (without data assimilation) and data assimilated SMS. For NLS, we take a reduced-order modeling approach, using one Gaussian mode to approximate the dynamics of the PDE. We incorporate measurements of the envelope of the solution of NLS into the DA-SMS solution. For KS and AD equations, we use shallow neural networks as our approximate solution and take pointwise measurements of the state of the system for DA-SMS. In each case, we add 5\% noise to the data to mimic real-world application where observational data is noisy.

\subsection{Nonlinear Schr\"{o}dinger Equation} \label{sec:NLS}
For our first example, we consider the nonlinear Schr\"{o}dinger equation. NLS models unidirectional, slowly modulating surface waves in deep fluids ~\cite{Anderson2022b,cousins16,Farazmand2017,Zakharov68}. Denoting the surface elevation by $\tilde{\eta}(\tilde{x},\tilde{t})$, we consider perturbations of a sinusoidal carrier wave,
\begin{equation}
\tilde{\eta}(\tilde{x},\tilde{t}) = Re\left\{\tilde{u}\left(\tilde{x},\tilde{t}\right)\exp \left[\hat{i} k_0 \tilde{x} - \omega_0 \tilde{t} \right]\right\},
\end{equation}
where $\hat i=\sqrt{-1}$, $\tilde{u}\left(\tilde{x},\tilde{t}\right)\in \mathbb{C}$ is the wave envelope, $k_0$ is the wave number of the carrier wave, and $\omega_0$ is its angular frequency. In dimensional variables, NLS reads~\cite{Zakharov68}
\begin{equation} \label{eq:dim_NLS}
\pard{\tilde{u}}{\tilde{t}} = -\hat{i} \frac{w_0}{8k_0^2} \frac{\partial^2 \tilde{u}}{\partial \tilde{x}^2} - \hat{i} \frac{\omega_0 k_0^2}{2} |\tilde{u}|^2 \tilde{u}, \quad \tilde{u}(x,0) = \tilde{u}_0(x).
\end{equation}
We introduce the non-dimensional variables $x = 2 \sqrt{2} k_0 \tilde{x}$, $t = - \omega_0 \tilde{t}$, and $u = (k_0/\sqrt{2})\tilde{u}$ as in~\cite{Anderson2022a}, so that Eq.~\eqref{eq:dim_NLS} becomes
\begin{equation}\label{eq:NLS}
\pard{u}{t} = \hat{i} \pard{^2u}{x^2} + \hat{i} |u|^2 u, \quad u(x,0) = u_0(x).
\end{equation} 
The associated boundary conditions for NLS are that the solution vanishes as $|x|$ approaches infinity. 

We solve \eqref{eq:NLS} using a Fourier pseudo-spectral truncation in space. In addition, we use an exponential time-differencing scheme as described in Ref.~\cite{Cox2002} with a time step of $\Delta t = 0.025$ as in Refs.~\cite{Anderson2022b,cousins16}. We discretize the spatial domain $\Omega = [-L/2,L/2]$ with length $L = 256\sqrt{2}\pi$, using $2^{11}$ Fourier modes with periodic boundary conditions. Our solutions are localized around $x=0$ and the domain is large enough to mimic the infinite domain size in NLS.

For SMS, our approximate solution takes the form,
\begin{equation}
\hat{u}(x,\params(t)) = A(t) \exp\left[-\frac{x^2}{L^2(t)} + \hat{i} \frac{V(t)}{L(t)}x^2 + \hat{i} \varphi(t)\right],
\end{equation}
in accordance with Refs.~\cite{Anderson2022a,Anderson2022b,Ruban2015a}. The corresponding shape parameters are given by $\params(t) = \{A(t),L(t),V(t),\varphi(t)\}$,
which control the amplitude $A$, length scale $L$, speed $V$ and phase $\varphi$ of the SMS solution. We evolve these parameters according to the SMS equation~\eqref{eq:SME_sym} which reads
\begin{equation}\label{eq:sms_NLS}
\dot{A} = -\frac{2AV}{L}, \quad \dot{L} = 4V, \quad \dot{V} = \frac{4}{L^3} - \frac{A^2}{\sqrt{2}L }, \quad \dot{\varphi} = \frac{5A^2}{4\sqrt{2}} - \frac{2}{L^2}. 
\end{equation}
These equations are integrated in time using Matlab's fourth-order Runge--Kutta scheme \texttt{ode45}. At the initial time, the parameters are $A_0 = 0.2$, $L_0 = 20$, $V_0 = 0$, and $\varphi_0 = 0$. For DNS we pick $u_0(x,0) = \hat{u}(x,\params(0))$, so that SMS and DNS solutions coincide initially. 

\begin{figure}[h!t]
\centering
\includegraphics[width = \textwidth]{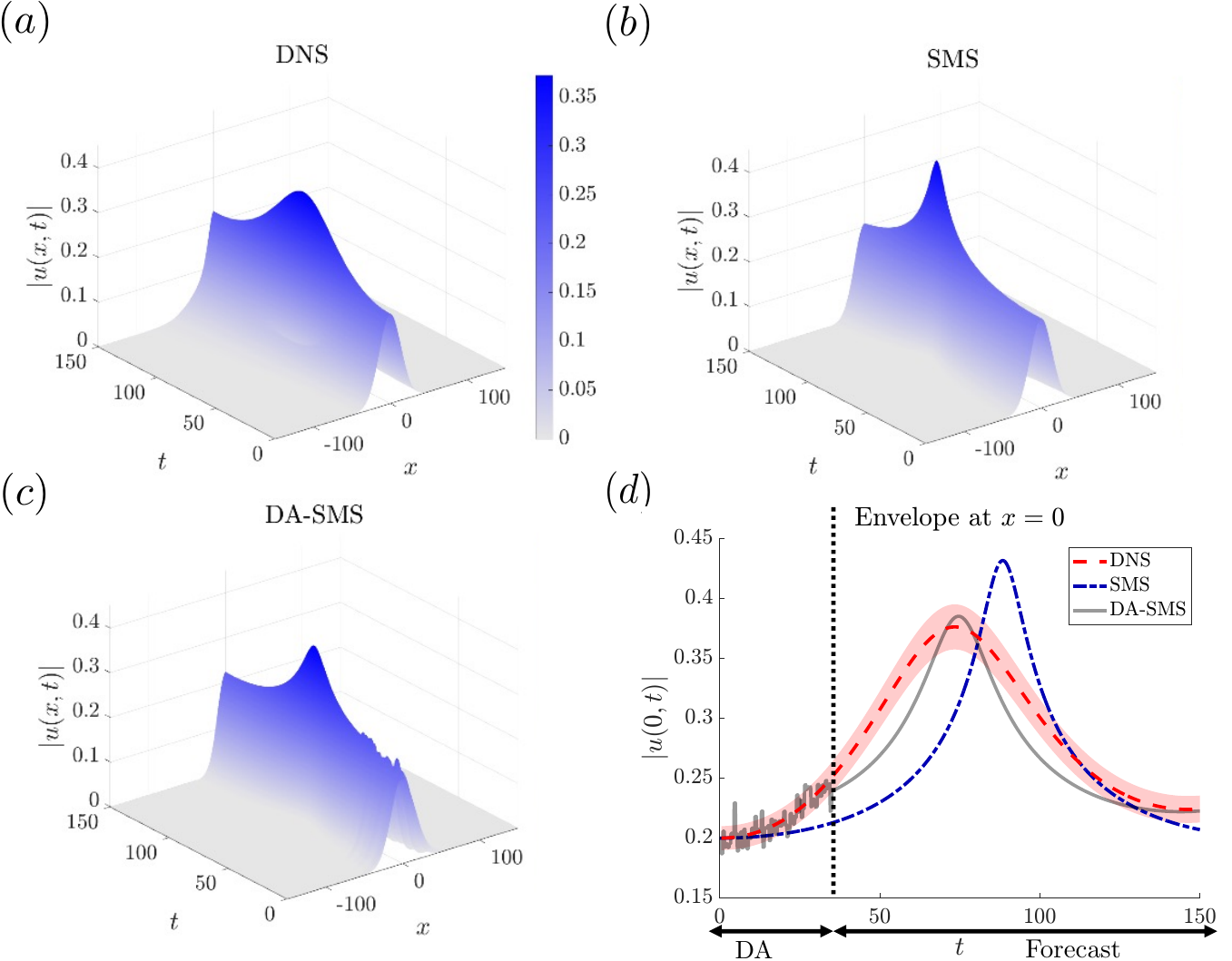}
\caption{\label{fig:NLS}Here we compare the modulus of the envelope $u(x,t)$ over the time interval $[0,150]$ using the initial parameters $A_0 = 0.2, \quad L_0 = 20, \quad V_0 = 0, \quad \phi_0 = 0$. We show (a.) DNS, (b.) SMS, (c.) DA-SMS. In (d.) we compare the envelopes of the three solutions at $x=0$, where DNS corresponds to the dashed red line, DA-SMS is the solid grey line, and SMS is the dot dashed blue line. The shaded red region denotes one standard deviation (5\%) away from the DNS solution.}
\end{figure}

For DA-SMS, at every $\Delta t = 0.5$ time units, we take pointwise measurements of the modulus of the envelope of the DNS solution at three sensor locations, $x_1=0$, $x_2 = 5$, and $x_3 = -10$. In other words, we have three pointwise observations $y_j(t_i) = |u(x_j,t_i)|$ at times $t_i\in\{0.5,1,1.5,\cdots\}$.
Since the solution is symmetric around $x=0$, the second and third sensor locations, $\{x_2,x_3\}$, are chosen asymmetrically to avoid redundant data.
The data-assimilation window is $\mathcal{T}_{DA} = [0,35]$, and the forecast window is $\mathcal{T}_f = [35,150]$.  We incorporate this data into the DA-SMS solution using Algorithm \ref{alg:DA-SMS}. 

Figure~\ref{fig:NLS} compares the solutions for DNS, DA-SMS, and SMS, as well as the  amplitudes of the solutions at $x=0$. For the chosen combination of initial amplitude $A_0$ and length scale $L_0$, we know that this wave exhibits a focusing behavior, meaning it increases in amplitude and decreases in width \cite{cousins16}. This focusing behavior is often used to model rogue waves on the ocean surface~\cite{chabchoub2011}.

As shown in Figure \ref{fig:NLS}(b,d), SMS correctly predicts that the wave will focus, but it overestimates the maximum amplitude. In addition, the maximum occurs with a time delay compared to the true DNS solution. In contrast, DA-SMS not only predicts the focusing of the wave, it also more accurately predicts the maximum amplitude. Moreover, DA-SMS correctly predicts the time at which the maximum wave height is attained; see Figure~\ref{fig:NLS}(c,d).
We note that the DA-SMS solution over the DA time window,  $\mathcal{T}_{DA} = [0,35]$, appears jagged because of the $5\%$ noise added to the observational data.

\subsection{Kuramoto--Sivashinsky Equation} \label{sec:KS}
For our second example, we consider the Kuramoto--Sivashinsky (KS) equation which models thermal instabilities in laminar flame fronts \cite{Kuramoto1978,Sivashinsky1977}. Here we approximate the solution to KS using a shallow neural network as the shape-morphing solution and compare the results to the DNS solution obtained from a Fourier pseudo-spectral method. KS equation is a fourth-order PDE,
\begin{equation}\label{eq:KS}
\pard{u}{t} = -u \pard{u}{x} - \pard{^2u}{x^2} - \pard{^4 u}{x^4}, \quad u(x,0) = u_0(x),
\end{equation}
with periodic boundary conditions on the domain $\Omega = [-L/2, L/2]$. Under certain conditions KS is known to exhibit chaotic behavior \cite{Kuramoto1978,Sivashinsky1977}. For our study, we choose $L = 22$ and define the function,
\begin{equation}\label{eq:KS_IC}
u^*_0(x) = \sin\left(\frac{2\pi}{L}x\right)+\sum_{k=2}^{4} \sin\left(\frac{2k\pi}{L} x\right) + \cos\left(\frac{2k \pi}{L}x\right).
\end{equation}
To ensure a reasonable, initial maximum amplitude, we rescale $u_0^*$ and define the initial condition,
\begin{equation}
u_0(x) = \frac{u^*_0(x)}{\max_{x\in\Omega} |u^*_0(x)|}.
\end{equation}
For DNS, we use a Fourier pseudo-spectral truncation using $2^7$ Fourier modes. We integrate in time using Matlab's fourth-order Runge--Kutta scheme \texttt{ode45}. 

On the other hand, for SMS we use a shallow neural network with a hyperbolic-tangent activation function in keeping with Anderson and Farazmand~\cite{Anderson2024a}. To enforce periodic boundary conditions, we use the methodology outlined in Ref. \cite{Du2021}, where the shape-morphing solution is given by
\begin{equation}\label{eq:SMS_KS}
\begin{split}
\hat{u}(x,\params(t)) & = \sum_{i=1}^N a_i(t) \tanh\left( w_i (t) s_i(x,c_i(t)) + b_i(t)\right), \\
 s_i(x,c_i(t)) &= \sin\left(\frac{2\pi }{L} x + c_i(t)\right). 
\end{split}
\end{equation}
In this case, the shape parameters, $\params(t)$, are the $N$ amplitudes $a_i(t)$, weights $w_i(t)$, and biases $\{b_i(t), c_i(t)\}$, so that $\params = \{a_i(t),w_i(t),b_i(t),c_i(t)\}_{i=1}^N$. Since symbolic computation is not possible in this case, we use the collocation method outlined in Section \ref{sec:prelim} to evolve the parameters in time. We use $2^7$ equidistant collocation points and a Tikhonov regularization parameter $\SMSTikh = 10^{-3}$ to counter the stiffness of the metric tensor. For the time stepping of the SMS equation, we again use Matlab's \texttt{ode45}. To find the initial parameters $\params_0$, we solve the nonlinear least squares problem, 
\begin{equation}
\min_{\params_0 \in \R^{4N}} \|\hat{u}(\cdot, \params_0) - u_0\|_{\mathcal{L}^2}^2, 
\end{equation}
using midpoint quadrature to discretize the integral in the $\mathcal{L}^2$ norm. For our study, we chose N = 10 modes, and found initial parameters $\params_0$ which correspond to an initial relative error less than 0.1\%. 
\begin{figure}[t!]
	\centering
	\includegraphics[width = \textwidth]{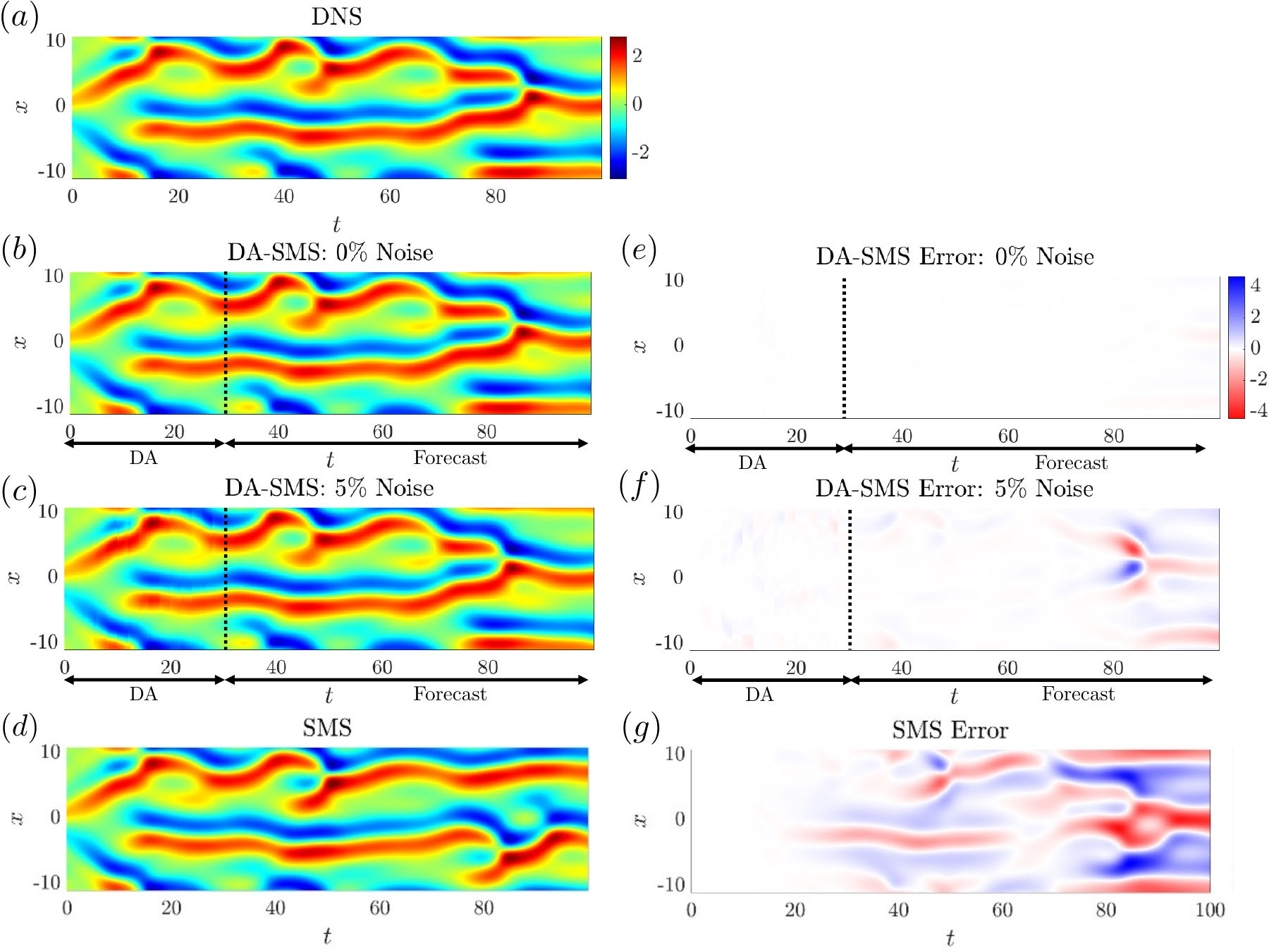}
	\caption{\label{fig:KS} Comparisons of the state of the system over the entire spatial domain $\Omega = [-11,11]$ and the entire duration of the simulation $t\in [0,100]$. The column on the left shows the results for: (a) DNS, (b) DA-SMS with clean (0\% noise) data, (c) DA-SMS with 5\% noise in the data, and (d) SMS with no data assimilation. In the column on the right we present the corresponding error defined by: $u_{DNS}(x,t) - \hat{u}(x,t)$ for (e) DA-SMS with clean data, (f) DA-SMS with noisy data, (g) SMS without data assimilation.}
\end{figure}

For DA-SMS, we take pointwise measurements of the state of the system at $\numsens = 10 $ locations as our observational data: $\data_j(t) = u(x_j,t)$ for $j= 1,2,\ldots, r$.  The sensors are equispaced throughout the domain, so that $x_j = -\frac{L}{2} + (j-1)\frac{L}{\numsens}$. We incorporate the data into the shape-morphing solution at every $\Delta t = 2.0$ time units, using the DA-SMS algorithm presented in Section \ref{sec:DA-SMS}. We assume we have data for the first 30 time units, so that $\mathcal{T}_{DA} = [0,30]$. We then forecast over the interval $\mathcal{T}_{f} = [30,100]$. For KS, we set the Tikhonov regularization parameter $\DATikh = 10^{-3}$ for DA, which is the same as the regularization parameter used in the collocation shape-morphing equation. We noticed that excellent agreement with observations are achieved with even a single Newton-like iteration. Therefore, to achieve computational speed up, we only take one Newton-like iteration in each step of Algorithm~\ref{alg:DA-SMS}. 

Figure \ref{fig:KS} compares DA-SMS results to both DNS and SMS without data assimilation. For DA-SMS, we consider two types of observational data: measurements without noise and noisy measurements where $5\%$ Gaussian noise is added to the observations $y_j$.
As shown in Figure \ref{fig:KS} (d,g), during the first 20 time units, the SMS solution closely matches the overall dynamics of the DNS solution. Soon after, however, the errors incurred by SMS move the trajectory of the approximate solution away from the true solution. More specifically, without data assimilation, the errors for SMS steadily grow from their minimum around $10^{-3}$, reach approximately one around $t=20$, and steadily grow for the remainder of the simulation.
This growth of errors over time is expected in a chaotic system. 

On the other hand, DA-SMS which incorporates data over the DA time window $t\in[0,30]$ keeps the approximate solution closer to the truth for a longer period of time. When there is no noise in the data, as shown in figure~\ref{fig:KS}(b,e), the errors over both DA time window and the forecast window remain near zero.
Even when 5\% error is added to observations (panels c and f), the error remains low for up to $t=75$. Eventually, as expected for a chaotic system, those errors also grow over the forecast time window where no data is available. However, compared to SMS without data assimilation, DA-SMS significantly extends the predictability horizon.
\begin{figure}[t!]
	\centering 
	\includegraphics[width = .9\textwidth]{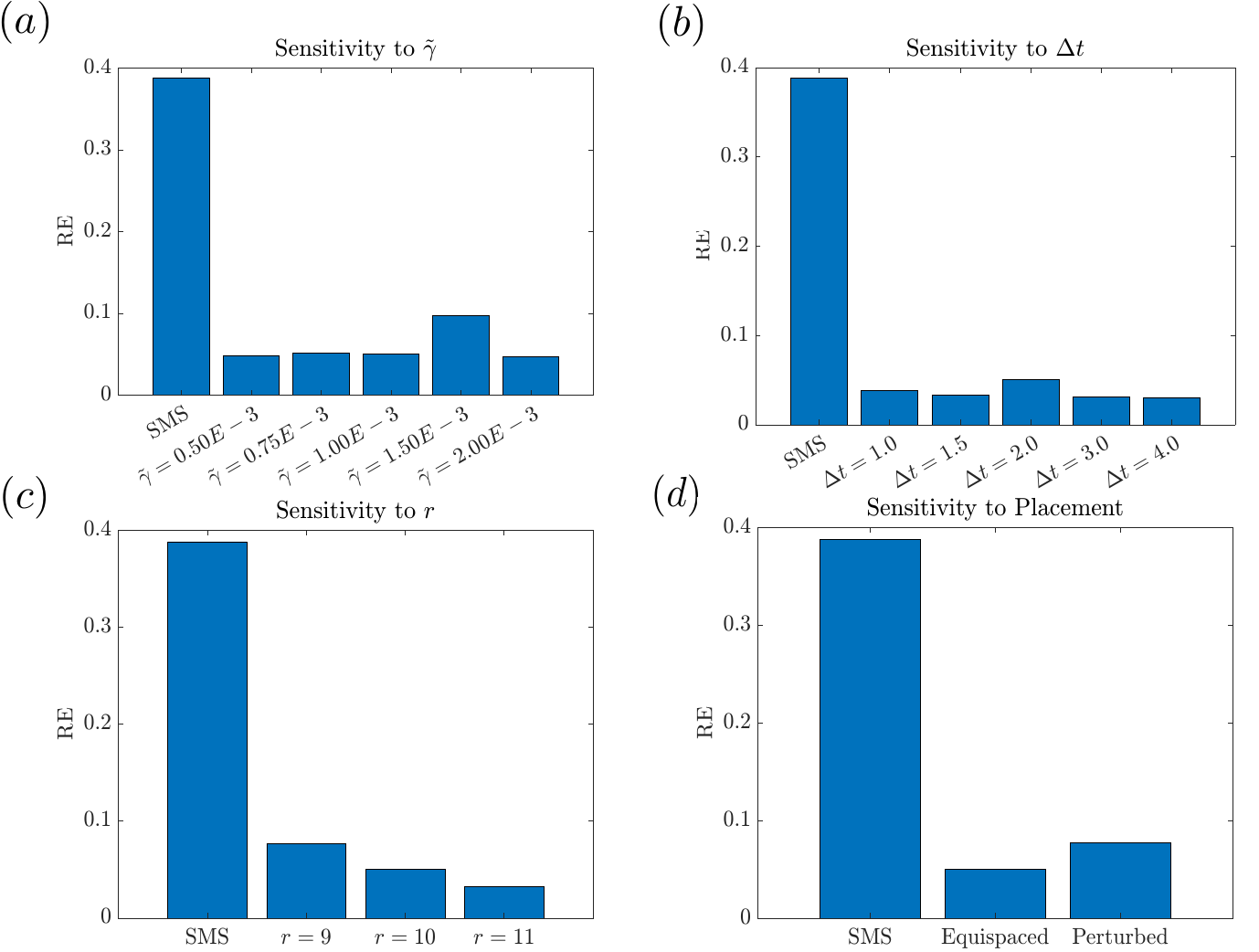}
	\caption{Comparison of the relative error (RE) for DA-SMS with different parameter values. The errors are reported at the end of the data assimilation time window ($t=30$). The first column in each panel shows the SMS error, without data assimilation, for comparison. (a) Varying the Tikhonov regularization parameter $\tilde{\gamma}$. (b) Different sampling rate $\Delta t$. (c) Different number of sensors $r$. (d) Perturbing the sensor locations.}
	\label{fig:sensitivity}
\end{figure}

Figure~\ref{fig:sensitivity} shows that our results are not too sensitive to the choice of parameters, such as the Tikhonov regularization parameter $\tilde{\gamma}$, the sampling rate $\Delta t$, number of observations $r$, or the placement of sensors. For instance, figure~\ref{fig:sensitivity}(a) shows the relative error of  DA-SMS at time $t=30$ (i.e., at the end of the data assimilation time window) for several values of the Tikhonov regularization parameter varying between $0.5\times 10^{-3}$ and $2\times 10^{-3}$. The relative errors are comparable for different values of the parameter, and significantly smaller than the SMS solution without data assimilation.
Similar results are observed in figure~\ref{fig:sensitivity}(b) where the time interval $\Delta t$ between observations varies.
Figure~\ref{fig:sensitivity}(c) shows that the relative error decreases as the number of sensors increases, as one would expect. 
Finally, figure~\ref{fig:sensitivity}(d) shows that equispaced sensors perform better. When perturbing the sensor locations, by randomly moving them one grid point to the left or to the right, the error increases slightly.

\subsection{Advection-Diffusion Equation} \label{sec:ADE}
\begin{figure}
\centering 
\includegraphics[width = \textwidth]{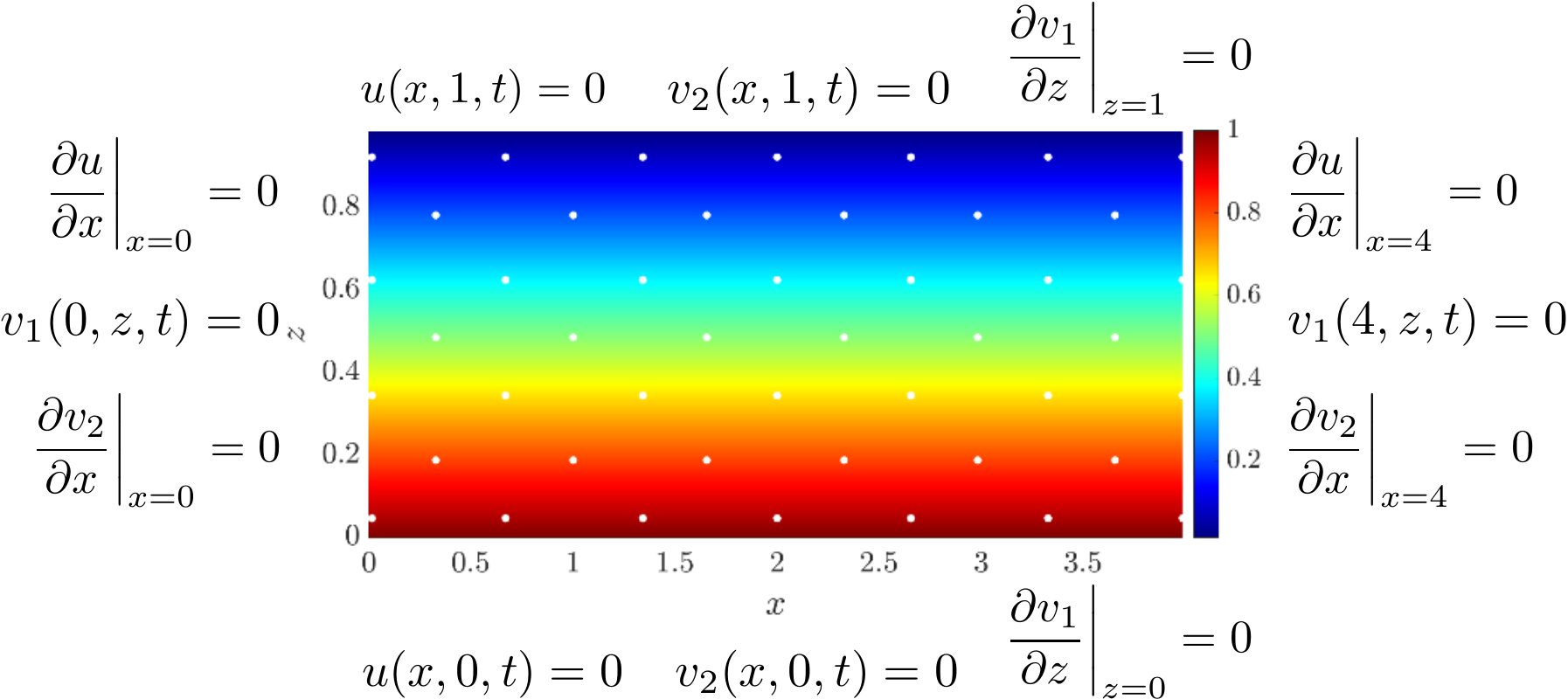}
\caption{Schematic illustration of the boundary conditions for the advection-diffusion equation in dimensionless variables. The domain size is $[0,4]\times [0,1]$.
The white dots indicate the sensor locations for DA-SMS.}
\label{fig:ADE_BCs}
\end{figure}

For our final example, we consider the passive advection-diffusion of temperature in a time-dependent fluid flow. Consider the temperature $\tilde{T}(\tilde{x},\tilde{z},\tilde{t})$ of the fluid in a rectangular domain $\Omega = [0,L]\times [0,H]$; see figure~\ref{fig:ADE_BCs}. At the horizontal boundaries, $z=0$ and $z=H$, the temperature is kept constant, so that $\tilde{T}(\tilde{x},0,\tilde{t}) = \tilde{T}_{bottom}$ and $\tilde{T}(\tilde{x},H,\tilde{t}) = \tilde{T}_{top}$, where $\tilde{T}_{bottom}>\tilde T_{top}$. For the left and right sides of the domain, the box is perfectly insulated so that no heat can escape in the $x$ direction, leading to the boundary conditions $\partial_{\tilde{x}}\tilde{T}(0,\tilde{z},\tilde{t}) = 0 = \partial_{\tilde{x}}\tilde{T}(L,\tilde{z},\tilde{t})$. The fluid is driven by a prescribed velocity field $\tilde{\ADEvel}(\tilde{x},\tilde{z},\tilde{t})$. We denote the horizontal and vertical components of the velocity by $\tilde{v}_1\left(\tilde{x},\tilde{z},\tilde{t}\right)$ and $\tilde{v}_2\left(\tilde{x},\tilde{z},\tilde{t}\right)$, respectively. 

In dimensional variables, the temperature evolves according to the advection-diffusion equation,
\begin{equation}\label{eq:AD_Ttil}
\pard{\tilde{T}}{\tilde{t}} + \tilde{\ADEvel} \cdot \nabla \tilde{T} = \tilde{\ADEconst} \Delta \tilde{T},
\end{equation}
where $\tilde{\ADEconst}$ is the thermal diffusivity of the fluid. If the fluid velocity is zero, $\tilde{\ADEvel}\equiv 0$, heat transfer occurs only by conduction, and the temperature profile of the system varies linearly with respect to $\tilde{z}$. More specifically, the conductive temperature profile $\tilde{T}_{c}$ is given by, 
\begin{equation} \label{eq:Tcond}
\tilde{T}_{c}\left(\tilde{z}\right) = \tilde{T}_{bottom}- \frac{\Delta \tilde{T}}{H} \tilde{z},
\end{equation}
where $\Delta \tilde{T} =  \tilde{T}_{bottom}-\tilde{T}_{top}$. If the fluid velocity is nonzero, the temperature will fluctuate away from the linear profile $\tilde{T}_c$. Thus it is convenient to instead study the evolution of the fluctuations $\tilde{\ADEvar}\left(\tilde{x},\tilde{z},\tilde{t}\right)=\tilde{T}\left(\tilde{x},\tilde{z},\tilde{t}\right) - \tilde{T}_{c}\left(\tilde{z}\right) $. Substituting this expression into \eqref{eq:AD_Ttil}, we obtain the following PDE for the fluctuations $\tilde{u}(\tilde{x},\tilde{z},\tilde{t})$,
\begin{equation}
\pard{\tilde{\ADEvar}}{\tilde{t}} + \tilde{\ADEvel} \cdot \nabla \tilde{\ADEvar} = \frac{\Delta \tilde{T}}{H} \tilde{v}_2 + \tilde{\ADEconst} \Delta \tilde{\ADEvar}. 
\end{equation}

We introduce the dimensionless variables,
\begin{equation}
 \vc x = \frac{\tilde{\vc x}}{H},\quad  T = \frac{\tilde{T}}{\Delta \tilde{T}}, \quad t = \frac{\tilde{t} v_f}{H}, \quad \text{and} \quad \tilde{\ADEvel} = \frac{\ADEvel}{v_f},
\end{equation}
where $\vc x = (x,z)$, and $v_f$ is the free-fall velocity given by $v_f = \sqrt{\beta \text{g} \Delta \tilde{T} H}$ where $\beta$ and $\text{g}$ are the thermal expansion of the fluid and the gravitational acceleration, respectively. In these dimensionless variables, the governing equation for the temperature fluctuation is given by
\begin{equation}\label{eq:ADE}
\pard{\ADEvar}{t} + \ADEvel\cdot \nabla \ADEvar = \ADEvelz + \ADEconst \Delta\ADEvar,
\end{equation}
where $\ADEconst = \tilde{\ADEconst} / \left(Hv_f\right)$. The corresponding boundary conditions for \eqref{eq:ADE} are homogeneous Neumann in $x$ and homogeneous Dirichlet in $z$, as shown in figure~\ref{fig:ADE_BCs}. 

Next, we specify the fluid velocity field. In order to be physically consistent, we must ensure that the velocity satisfies the appropriate boundary conditions. Namely, we require that there are no slip conditions on the horizontal component, $\ADEvelx$, at the left and right boundaries and no-flux conditions in the $z$ direction on the top and bottom. Additionally, the vertical component of the velocity $\ADEvelz$ should have no slip conditions on the top and bottom of the domain, and no-flux conditions in the $x$ direction on the left and right of the domain. In summary, the fluid velocity field must satisfy the following boundary conditions,
\begin{equation}\label{eq:v_bc}
\begin{split}
\ADEvelx(x,z,t)\bigg|_{x = 0, L} &= 0, \quad \pard{\ADEvelx}{z}\bigg|_{z = 0, H} = 0, \\
 \ADEvelz(x,z,t)\bigg|_{z = 0,H} &= 0, \quad \pard{\ADEvelz}{x}\bigg|_{x = 0, L} = 0. 
\end{split}
\end{equation}
We summarize the boundary conditions for the temperature fluctuations $\ADEvar$ and the fluid velocity $\ADEvel$ in figure~\ref{fig:ADE_BCs}. 
Furthermore, we require $\nabla \cdot \ADEvel = 0$, ensuring that the fluid flow is incompressible.

To attain these requirements on $\ADEvel$, we define the velocity field,
\begin{equation}
\ADEvel(x,z,t) =\pi A \begin{pmatrix}
-\sin(\pi f(x,t)) \cos(\pi z) \\
\cos(\pi f(x,t))\sin(\pi z) \pard{f}{x}
\end{pmatrix}
\end{equation}
where 
\begin{equation}
f(x,t) = m \frac{x}{L} + \epsilon L^4 \sin(\omega t) \left[\frac{x}{L} - 2 \left(\frac{x}{L}\right)^3 + \left(\frac{x}{L}\right)^4\right].
\end{equation}
This velocity field, which satisfies the boundary conditions~\eqref{eq:v_bc},  is a modified version of double gyre flow that is often used in oceanography~\cite{computeVariLCS,shadden05}.  
The parameter $m$ roughly determines the number of gyres (or vortices) in the $x$ direction. The parameter $\omega$ is the frequency of gyre oscillations
and the small parameter $\epsilon$ determines the amplitude of these oscillations. For our simulations we set $m=2$, $A = 0.1$, $\epsilon = 0.025$, and $\omega = \pi$.  For the spatial domain, we set the dimensionless width of the domain $L = 4$ and the height $H = 1$. We set the temperatures at the bottom and top of the domain as $T_{bottom} = 1$ and $T_{top} = 0$, respectively, and we set $\ADEconst = 10^{-3}$ for the dimensionless thermal diffusivity. Lastly, we choose the initial condition, 
$$
\ADEvar_0(x,z) = 10^{-1}\cos\left(\frac{\pi}{L} x \right)\sin\left(\frac{\pi}{H}z\right).
$$
\begin{figure}[ht]
	\centering
	\includegraphics[width = \textwidth]{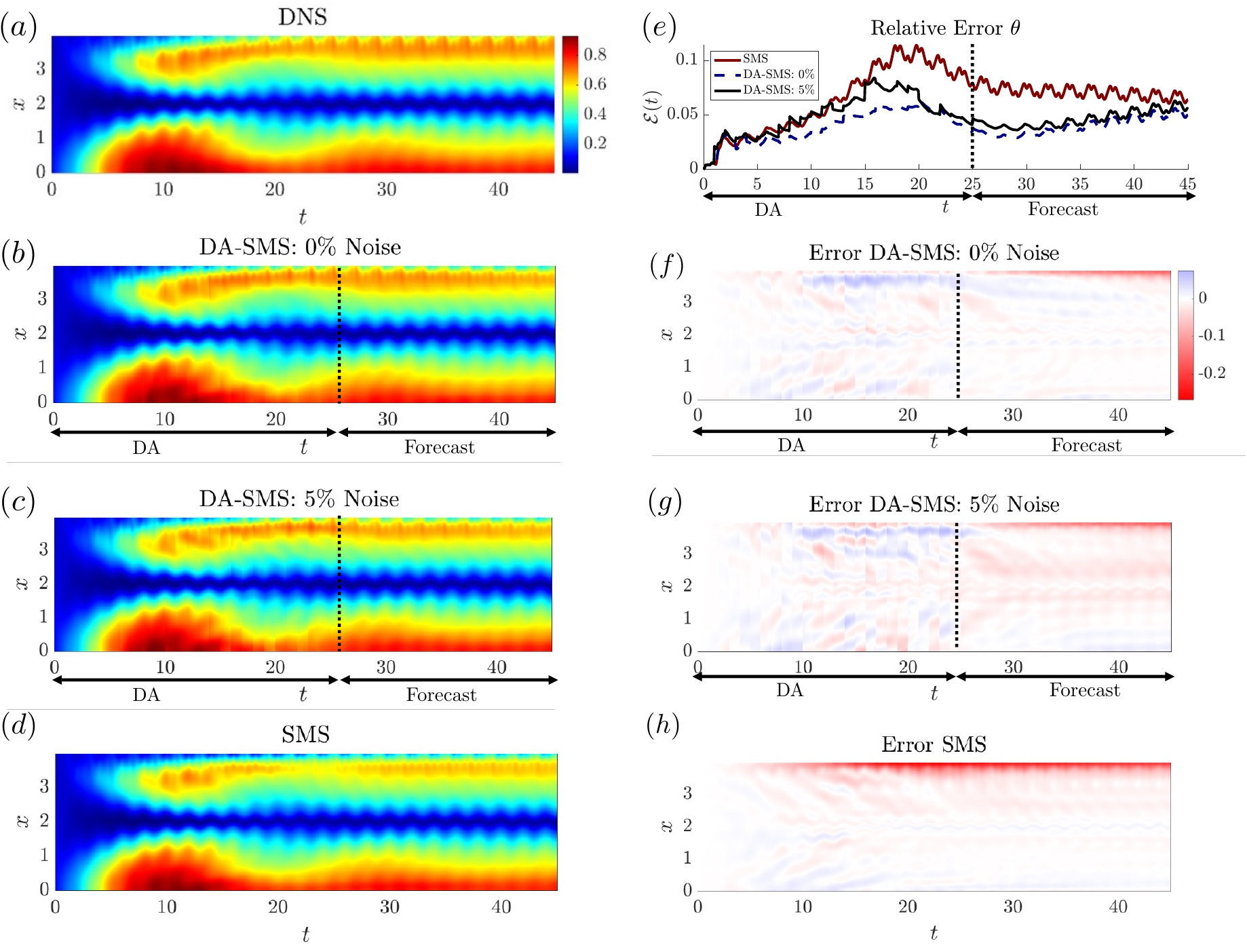}
	\caption{\label{fig:ADE} Space-time plots of the solutions to ADE \eqref{eq:ADE} where we have taken one slice in space located at $z = 0.9$ at every time instance over the interval $[0,45].$ The left column shows the results obtained from (a) DNS, (b) DA-SMS with $0\%$ noise in the data, (c) DA-SMS with $5\%$ noise in the data, and (d) SMS without data assimilation. Panels (f,g,h) show the corresponding difference $u_{DNS}(x,0.9,t) - \hat{u}(x,0.9,\params(t))$. Panel (e) shows the total relative errors.}
\end{figure}

For DNS, we use a Fourier pseudo-spectral method to solve equation~\eqref{eq:ADE}, using $2^8$ modes in the $x$ direction, and $2^6$ modes in the $z$ direction. The DNS solution is shown in figure~\ref{fig:ADE}(a). In order to visualize the space-time evolution of the solution, we only show the horizontal slice at $z=0.9$.

For the shape-morphing solution $\hat u$, we use a shallow neural network of width $N=100$ as described in~\ref{app:2DNN}. In solving PDEs with  neural networks, the boundary conditions are often approximately enforced by adding them to the loss function as a penalty term~\cite{Raissi2019}. Here, in contrast, we construct our neural network such that the boundary conditions are explicitly built into the network and therefore are exactly satisfied for any set of network parameters $\params$. The details of our construction can be found in~\ref{app:2DNN}.

The evolution of network parameters $\params(t)$ is determined by solving the SMS equation~\eqref{eq:reg_coll} with Tikhonov regularization parameter $\SMSTikh = 5\times 10^{-2}$ and using $2^6 \times 2^4$ collocation points evenly spaced throughout the domain $[0,4]\times [0,1]$. To initialize the SMS equation, we solve the nonlinear optimization problem,
$\min_{\params_0\in \mathbb{R}^{6N}}\|\hat{\ADEvar}(\cdot,\params_0) - \ADEvar_0\|_{\mathcal{L}^2}^2$,
to obtain the initial condition $\params_0$. The resulting approximate initial condition $\hat u(\cdot,\params_0)$ matches the true initial condition $u_0$ within $0.1\%$ error. Lastly, the temporal integration of the SMS equation is carried out using \texttt{ode45}. The resulting SMS solution without data assimilation is shown in figure~\ref{fig:ADE} (d,h).
 
Now we turn our attention to DA-SMS, where the observables are pointwise measurements $y_i(t) = \ADEvar(x_i,z_i,t)$. We place $\numsens = 46$ sensors on a scattered grid as shown in figure~\ref{fig:ADE_BCs}. The observations are gathered every $\Delta t = 0.5$ time units over the data assimilation window $\mathcal{T}_{DA}  = [0,25]$. The forecast window is $\mathcal{T}_{f} = [25,45]$.  For the Newton iterations, we set the Tikhonov regularization parameter to $\DATikh = 5\times 10^{-2}$, which is equal to the Tikhonov regularization parameter for SMS. Similar to the KS equation, we observed that even one Newton iteration suffices to obtain satisfactory agreement with observations.

Figure~\ref{fig:ADE} compares the results of DA-SMS and SMS without data assimilation . Panel (e) shows the relative $\mathcal L^2$ error, 
\begin{equation}\label{eq:RE}
	\mathcal{E}(t) = \frac{\|u(\cdot, t) - \hat{u}(\cdot, \params(t))\|_{\mathcal L^2}}{\|u(\cdot,t)\|_{\mathcal L^2}}, 
\end{equation}
for each simulation. Although the SMS solution captures the overall dynamics, its relative error increases to about 12\% over time. Recall that the boundary conditions are exactly enforced in the construction of our neural network, so that the error is exactly zero at the boundary. However, immediately near the boundary we observe a large accumulation of error for the SMS solution. This is particularly pronounced near $x=4$. On the other hand, the DA-SMS is able to suppress these errors significantly. For DA-SMS with clean data (no noise added to the observations), the total relative error is at most 5\%. When noise is added to the observational data, the DA-SMS error increases slightly to at most 7\%. However, this error is still smaller than the SMS simulations without data assimilation.

\section{Conclusions}\label{sec:conc}
Shape-morphing solutions are a computational method for obtaining numerical solutions to PDEs. In spite of the rapidly growing literature on this subject, data assimilation with SMS had remained unaddressed. Here, we developed two data assimilation methods for shape-morphing solutions: (i) Discrete-time data assimilated SMS and (ii) Continuous-time data assimilated SMS. 
Both methods seek to use observations to correct for modeling error, approximation error, and uncertainties in the initial condition.

Our discrete method takes the form of a predictor-corrector scheme. During the prediction phase, the SMS equations are used to evolve the solution between consecutive observation points. At time instances where observations are available, a Newton-like iteration is used to correct for the discrepancy between the SMS and the observational data. We proved that, under certain conditions, the DA-SMS method converges uniformly to the true solution.

Our continuous-time method takes a different approach. It treats the observational data as constraints to the optimization problem underlying SMS equations. The resulting constrained SMS equation ensures that the solution agrees with the observational data up to a constant which is independent of time.

Several open questions remain to be addressed. These include:

(i) Better convergence analysis: Our theoretical results concerning the convergence of the discrete-time DA-SMS are rather pessimistic. They indicate that a dense set of observational data is required for the method to converge. In contrast, our numerical results indicate that the method works with relatively few observations and a single Newton-like iteration. Thus, our numerical results suggest that tighter error bounds for the DA-SMS error are feasible.

(ii) Sensor placement: It is well-known that the quality of data assimilation methods depend significantly on the location of the sensors which gather the data~\cite{Attia2018,farazmand2024}.
Here, we placed our sensors in an ad hoc manner. Specialized sensor placement methods, tailored for use with shape-morphing solutions, are desirable.

(iii) Other data assimilation techniques: Classical data assimilation methods, such as extended Kalman filtering, variational data assimilation, and nudging techniques can in principle be used with SMS. However, the nonlinear dependence of SMS on its parameters $\params(t)$ complicates a straightforward implementation of these methods. Future work should explore the necessary modifications to existing DA methods for use with SMS. Furthermore, the nonlinearities necessitate a separate convergence analysis for existing DA methods used in conjunction with SMS.

\subsubsection*{Acknowledgments}
This work was supported by the National Science Foundation under the grants DMS-1745654 (Research Training Group), DMS-2208541 (Computational Mathematics Program), and DMS-2220548 (Algorithms for Threat Detection). 

\appendix
\section{Derivation of Continuous-time Data Assimilation Equations}\label{app:cont_deriv}
Here we present a detailed derivation of the results presented in Section \ref{sec:cont_DA} for both types of shape-morphing equation. To begin, we consider the optimization problem
\begin{align*}
\min_{\dot{\params}} & \quad \frac{1}{2}\hilnorm{R(\vc x, \params, \dot{\params})}^2 + \frac{\SMSTikh}{2}\|\dot{\params}\|_2^2 \\ 
\text{s.t.} & \quad \jacop \dot{\params} = \dot{\vc \data}. 
\end{align*}
To solve the constrained optimization problem we introduce the Lagrange multipliers $\vc \lambda \in \R^r$. To this end, we construct the constrained cost functional 
\begin{equation}
\mathcal{J}_c(\params; \dot{\params}, \vc \lambda) = \frac{1}{2}\hilnorm{R(\vc x, \params, \dot{\params})}^2  + \vc \lambda^\top (\jacop\dot{\params} - \dot{\vc \data}) + \frac{\SMSTikh}{2}\|\dot{\params}\\|_2^2
\end{equation}
Next, we expand the Hilbert norm in terms of the Hilbert inner product, which yields
\begin{equation}\label{eq:ip_opt}
\begin{split}
\mathcal{J}_c(\params; \dot{\params}, \vc \lambda) = & \frac{1}{2} \hilip{\sum_{i=1}^{\numparams N} \pard{\hat{u}}{\singparam_i}\dot{\singparam_i}}{\sum_{j=1}^{\numparams N} \pard{\hat{u}}{\singparam_j}\dot{\singparam_j}} - \hilip{\sum_{i=1}^{\numparams  N} \pard{\hat{u}}{\singparam_i}\dot{\singparam_i}}{F(\hat{u})} \\
& - \frac{1}{2} \hilnorm{F(\hat{u})}^2 + \vc \lambda^\top (\jacop\dot{\params} - \dot{\vc \data}) + \frac{\SMSTikh}{2}\|\dot{\params}\|_2^2. 
\end{split}
\end{equation}
Following Ref. \cite{Anderson2022a}, we introduce the metric tensor $M$ and the vector $\vc f$ which we describe in Section \ref{sec:prelim}, so we can re-express \eqref{eq:ip_opt} as 
\begin{equation}
\mathcal{J}_c(\params; \dot{\params}, \vc \lambda) = \frac{1}{2}\dot{\params}^\top \left(M + \gamma I\right )\dot{\params} -\dot{\params}^\top \vc f + \hilnorm{F(\hat{u})}^2 + \vc \lambda^\top (\jacop\dot{\params} - \dot{\vc \data}). 
\end{equation}
Next, we compute the gradient of $\mathcal J_c$ with respect to $\dot{\params}$ and set it equal to zero to solve for $\dot{\params}$,
\begin{equation}
\nabla_{\dot{\params}} \mathcal{J}_c = (M+\SMSTikh I)\dot{\params} - \vc f + \jacop^\top \vc \lambda=\vc 0 
\end{equation}
By construction $M$ is symmetric positive semidefinite; because $\SMSTikh >0$, $(M+\SMSTikh I)$ is invertible, thus 
\begin{equation}\label{eq:app_thetadot}
\dot{\params} = \left(M+\SMSTikh I \right)^{-1}\left(\vc f - \jacop ^\top \vc \lambda\right). 
\end{equation}
Next, we take the gradient of $\mathcal J_c$ with respect to $\vc \lambda$ and substitute the expression found in \eqref{eq:app_thetadot} to solve for $\vc \lambda$. Letting $M_\SMSTikh = M + \SMSTikh I$ we have, 
\begin{equation}
\jacop \left(M_\SMSTikh ^{-1}\vc f - M_\SMSTikh ^{-1}\jacop^\top\vc \lambda \right) = \dot{\vc \data}
\end{equation}
Isolating $\vc \lambda$, we find 
\begin{equation}
\jacop M_\SMSTikh ^{-1}\jacop^\top\vc \lambda  =\jacop M_\SMSTikh^{-1} \vc f-  \dot{\vc \data}
\end{equation}
In most cases $\numsens \ll kN$, so $\jacop$ is a short-fat matrix. Thus if $\jacop$ is full rank, then $\jacop M_\SMSTikh^{-1} \jacop^\top$ is symmetric positive definite. Hence,
\begin{equation}
\vc \lambda = \left(\jacop M_\SMSTikh^{-1} \jacop^\top\right)^{-1}\jacop M_\SMSTikh^{-1} \vc f-  \left(\jacop M_\SMSTikh^{-1} \jacop^\top\right)^{-1} \dot{\vc \data}
\end{equation}
Lastly, we use this expression for $\vc \lambda$ in Eq. \eqref{eq:app_thetadot}, and we arrive at 
\begin{equation}
\dot{\params}  = \left[I - M_\SMSTikh^{-1}\jacop^\top \left(\jacop M_\SMSTikh^{-1} \jacop^\top \right)^{-1}\jacop \right]M_\SMSTikh^{-1} \vc f + M_{\gamma}^{-1}\jacop^\top\left(\jacop M_\gamma^{-1}\jacop^\top\right)^{-1}\dot{\vc \data} . 
\end{equation}

Next, we turn our attention to the collocation shape-morphing equation. To begin we amend the collocation optimization problem in Section \ref{sec:prelim} by adding the same constraint as above. In this case we have 
\begin{equation}
\begin{split}
\min_{\dot{\params\in \R^{Nk}}} & \  \frac{1}{2} \|\tilde{M}\dot{\params} - \tilde{\vc f}\|^2_2 + \frac{\gamma}{2}\|\dot{\params}\|_2^2, \\
\text{s.t.} & \ \jacop \dot{\params} = \dot{\vc y}
\end{split}
\end{equation}
To find the solution we follow the same procedure as above, first introducing the Lagrange multipliers $\vc\lambda\in\R^\numsens$ and constructing the constrained cost function, 
\begin{equation}
\tilde{\mathcal{J}}_c(\params; \dot{\params}, \vc \lambda) =  \frac{1}{2} \|\tilde{M}\dot{\params} - \tilde{\vc f}\|_2^2 + \frac{\gamma}{2}\|\dot{\params}\|_2^2 + \vc \lambda^\top \left(\jacop \dot{\params} - \dot{\vc y}\right). 
\end{equation}
Next we expand the cost function, compute the derivative with respect to $\dot{\params}$, 
\begin{equation}
\nabla_{\dot{\params}} \tilde{\mathcal{J}}_c = \left(\tilde{M}^\top \tilde{M}+\SMSTikh I\right) \dot{\params} - \tilde{M}^\top\tilde{\vc f} + \jacop^\top \vc \lambda 
\end{equation}
To find the critical point, we set this expression to zero and solve for $\dot{\params}$. Since $\gamma>0$ we have that $\tilde{M}_\SMSTikh = \tilde{M}^\top\tilde{M}+\SMSTikh I$ is invertible, so 
\begin{equation}
\dot{\params} = \tilde{M}_\SMSTikh^{-1} \tilde{M}^\top \tilde{\vc f} - \tilde{M}_\SMSTikh^{-1} \jacop^\top\vc \lambda
\end{equation}
Next, we compute the gradient with respect to $\vc \lambda$, which recovers the constraint, and then substitute our expression for $\dot{\params}$. If $\jacop$ is full rank, then $\jacop \tilde{M}_\gamma^{-1}\jacop$ is invertible, and we find that 
\begin{equation}
\vc \lambda = \left(\jacop \tilde{M}_{\SMSTikh}^{-1} \jacop \right)^{-1} \left(\jacop \tilde{M}_\SMSTikh^{-1} \tilde{M}^\top \tilde{\vc f}-\dot{\vc \data}\right). 
\end{equation}
Combining these expressions and simplifying, we find that the continuous data assimilation algorithm for the collocation method is given by 
\begin{equation}\begin{split}
\dot{\params} = & \left[I -\tilde{M}_{\SMSTikh}^{-1}\jacop^\top \left(\jacop \tilde{M}_{\SMSTikh}^{-1}\jacop^\top\right)^{-1}  \right] \tilde{M}_{\SMSTikh}^{-1}\tilde{M}^\top\tilde{\vc f} \\
& + \tilde{M}_\SMSTikh^{-1}\jacop^\top \left(\jacop \tilde{M}_{\SMSTikh}^{-1}\jacop^\top\right)^{-1} \dot{\vc \data}
\end{split}
\end{equation}

\section{Neural Network Construction for Mixed Dirichlet-Neumann Boundary Conditions}\label{app:2DNN}
In this section we discus the construction of the neural network used for the advection diffusion equation in Section~\ref{sec:ADE}. This neural network is designed such that the mixed Dirichlet-Neumann boundary conditions of the AD equation are automatically satisfied by the network. This objective is achieved in three steps:
\begin{enumerate}
	\item First, we define an extended computational domain $\Omega_e$ using copies of the physical rectangular domain $\Omega$.
	\item We then enforce periodic boundary conditions on the extended domain $\Omega_e$.
	\item Finally, we enforce odd and even symmetries into the network to ensure the mixed Dirichlet-Neumann boundary conditions are satisfied on the physical domain $\Omega$.
\end{enumerate}

Consider a general neural network $N(\vc x, \params)$, which may be deep or shallow. The neural network is defined on the physical spatial domain $\vc x\in\Omega = [0,L]\times [0,H]$. Recall that the AD equation has homogeneous Dirichlet boundary conditions in $z$ and homogeneous Neumann boundary conditions in $x$. First, we define the extended computational domain $\Omega_e = [-L,L]\times[-H,H]$. We impose periodic boundary conditions over the entire computational domain $\Omega_e$ using the methodology presented in Ref.~\cite{Du2021}. Namely, in each node of the network $N(\vc x, \params)$, we apply a change of variables through the function,
\begin{equation}
\vc s_i(\vc x) = \left[\sin\left(\frac{\pi }{L} x + c^x_i \right),\sin\left(\frac{\pi }{H} z + c^z_i\right) \right]^\top,
\end{equation}
where $c^x_i$ and $c^z_i$ are translation parameters.

For example, in AD equation where we consider a shallow neural network with a hyperbolic tangent activation function, we have
\begin{equation}
N_p(\vc x, \params) =\sum_{i=1}^N a_i\tanh(\vc w_i \cdot \vc s_i(\vc x) + b_i), 
\end{equation}
where $\vc w_i = (w_i^x,w_i^z)$ are network weights and $b_i$ are its biases.
This neural network has periodic boundary conditions on the extended domain $\Omega_e$. The complete set of network parameters is therefore given by $\params = \{a_i, b_i, w_i^x, w_i^z, c_i^x, c_i^z\}_{i=1}^N$.

To enforce the mixed Dirichlet-Neumann  boundary conditions on the physical domain $\Omega$, we exploit the inherent symmetries of even and odd functions. For a continuous odd function, $f(-x) = -f(x)$, we have $f(0)=0$. Additionally, for a continuously differentiable even function, $f(-x) = f(x)$, its derivative must be odd, and therefore, $f'(0)=0$. As a result, if the construction of the neural network is odd with respect to $z$ and even with respect to $x$, then it satisfies the correct boundary conditions of the AD equation. 

Finally, we introduce these symmetries by using the neural network $N_p(\vc x,\params)$ to define the SMS solution,
\begin{equation}
\hat{u}(\vc x, \params) = N_p(x,z,\params) - N_p(x,-z,\params) + N_p(-x,z,\params) - N_p(-x,-z,\params). 
\end{equation}
It is straightforward to verify that $\hat{u}$ is odd with respect to $z$ and even with respect to $x$. This, together with periodicity of $\hat u$ on the extended computational domain $\Omega_e$, ensures that the SMS solution $\hat u$ satisfies the mixed homogeneous Dirichlet-Neumann boundary conditions on the physical domain $\Omega$.


\end{document}